\documentclass[11pt]{article}
\usepackage{amsmath,amsthm,amssymb,amscd,fancybox,ifthen,float,epsfig,subcaption}
\usepackage[all]{xy}
\usepackage{times}
\usepackage{hyperref}

\floatstyle{plain}
\restylefloat{figure}

\usepackage{color}

\definecolor{gr}{rgb}   {0.,   0.8,   0. }
\definecolor{bl}{rgb}   {0.,   0.5,   1. }
\definecolor{mg}{rgb}   {0.7,  0.,    0.7}

\newcommand{\Bk}{\color{black}}
\newcommand{\Rd}{\color{black}}

\newcommand\mf{\operatorname{mf}}
\newcommand\Sc{\operatorname{sc}}

\newcommand\ff{\operatorname{ff}}

\newcommand\LAP{\operatorname{LAP}}
\newcommand\WG{\operatorname{WG}}

\newcommand\Diff{\operatorname{Diff}}
\newcommand\Th{\operatorname{th}}

\newcommand\WF{\operatorname{WF}}

\newcommand\tot{\operatorname{tot}}

\newcommand\In{\operatorname{in}}
\newcommand\Int{\operatorname{int}}

\newcommand\out{\operatorname{out}}

\newcommand\cJ{\mathcal J}

\newcommand\cV{\mathcal V}

\newcommand\cR{\mathcal R}

\newcommand\tx{\tilde x}

\newcommand\tm{\widetilde m}

\newcommand\tvarphi{\widetilde \varphi}

\newcommand\tu{\widetilde u}

\newcommand\cF{\mathcal F}
\newcommand\cE{\mathcal E}
\newcommand\cW{\mathcal W}

\newcommand\hN{\widehat{N}}

\newcommand\ka{\kappa}

\newcommand\PV{\operatorname{P.V.}}

\newcommand\cA{\mathcal{A}}
\newcommand\cC{\mathcal{C}}

\newcommand\cS{\mathcal{S}}

\newcommand\cD{\mathcal{D}}

\newcommand\hv{\hat v}

\newcommand\hf{\hat{f}}
\newcommand\hh{\hat{h}}

\renewcommand\Im{\operatorname{Im}}

\newcommand\bbC{\mathbb C}

\newcommand\bbN{\mathbb N}

\newcommand\bbR{\mathbb R}

\newcommand\pa{\partial}

\newcommand\restrictedto{\upharpoonright}

\newcommand\supp{\operatorname{supp}}
\newcommand\singsupp{\operatorname{singsupp}}
\newcommand\scsingsupp{\operatorname{sc-singsupp}}
\newcommand\subsubset{\subset\!\subset}

\newcommand\Id{\operatorname{Id}}

\newtheorem{theorem}{Theorem}

\newtheorem{lemma}{Lemma}

\theoremstyle{definition}
\newtheorem{definition}{Definition}

\newtheorem{example}{Example}

\theoremstyle{remark}
\newtheorem{remark}{Remark}

\begin{document}

\title{Solving the Scattering Problem for Open Wave-Guide Networks, III:
\\Radiation Conditions and Uniqueness}

\author{Charles L. Epstein\footnote{ Center for Computational Mathematics, 
  Flatiron Institute, 162 Fifth Avenue, New York, NY 10010. E-mail:
  {cepstein@flatironinstitute.org}.}\,\,   and Rafe Mazzeo\footnote{Dept. of
  Mathematics, Stanford University, Stanford, CA. E-mail: rmazzeo@stanford.edu}
  }
\date{November 6, 2025}

\maketitle
\begin{abstract} This paper continues the analysis of the scattering problem for 
  a network of open wave-guides started in~\cite{EpWG2023_1,EpWG2023_2}. In this
  part we present explicit, physically motivated radiation conditions that ensure uniqueness 
  of the solution to the scattering problem. These conditions stem from a 2000 paper of 
  Vasy on 3-body Schr\"odinger operators, see~\cite{VasyAsterisque}; we also discuss 
  closely related conditions from a 1994 paper of Isozaki ~\cite{Isozaki94}. Vasy's paper 
  also proves   the existence of the limiting absorption resolvents, and that the limiting
  solutions satisfy the radiation conditions.  The statements of these results require a
  calculus of pseudodifferential operators, called the 3-body scattering calculus, which
  is briefly introduced here. We show that the solutions to the model problems
  obtained in~\cite{EpWG2023_1} satisfy these radiation conditions, which makes it
  possible to prove uniqueness, and therefore existence, for the system of
  Fredholm integral equations introduced in that paper.
\end{abstract}
\tableofcontents

\section{Introduction}\label{sec1}
Many opto-electronic and photonic devices are modeled as \emph{open wave-guide
networks}.\footnote{\Rd In the Applied Math, Engineering and Physics literature
an ``open wave-guide'' usually refers to a translationally invariant device.  We
call these {\em bi-infinite wave-guides.} The main point of our work is that we
consider an assemblage of devices that are asymptotically modeled by bi-infinite
wave-guides, which we call a {\em wave-guide network.} \Bk} In such devices
there are channels\footnote{This usage of the term `channel' is standard for
wave-guides. It is different from the terminology used in the $N$-body
Schr\"odinger equation literature, where channels refer to distinguished
eigenspaces of subsystems. What we call channels are analogous to `collision
planes' in the Schr\"odinger equation literature.} defined by spatial variations
in the electrical permittivity, but the channels are unclad, so the
electromagnetic waves are not confined to the channel. Such physical systems are
described by Maxwell's equations with spatially dependent permittivity.

In this paper we consider a simpler scalar model, wherein the permittivity
$\epsilon(x)$ is a positive real-valued function, equal to the positive constant
$\epsilon_1$  outside the set $\Omega,$ described below. \Rd
{The underlying physical model is a wave equation
\begin{equation}\label{eqn3.201}
  \Delta U=\epsilon(x)\pa_t^2U,
\end{equation}
where we assume that the solution takes the form
$U(x,t)=e^{-ift}u(x),$ for an $f>0.$   That is, we  consider the time harmonic
case and seek solutions to
\begin{equation}\label{eqn3.200}
(\Delta+k^2(x))u=0,
\end{equation}
where
\begin{equation}
k^2(x)=\epsilon(x)f^2.
\end{equation}
\Rd The function $\epsilon(x)$ can also interpreted as the reciprocal of the square
of the ``sound speed.'' With this interpretation, equation~\eqref{eqn3.200} is
model for acoustic scattering from an acoustic wave-guide network.} \Bk

The set, $\Omega,$ where $\epsilon(x)$ differs from the free-space value,
$\epsilon_1,$ is the union of a compact set and finitely many `tubes' extending
to infinity. More precisely 
$$\Omega=\Omega_0\cup\left[\bigcup_{\alpha\in\cA}
T_{\alpha}\right],$$
where $\cA$ is a finite index set, and $\Omega_0$ is
a compact set. Each $T_{\alpha}$ is a tube, which is unbounded in one direction:
Fix a set of points $\{v_{\alpha}:\:\alpha\in\cA\}$ on the unit sphere
$S^{d-1}$, for some positive constants, $d_\alpha, R_\alpha$, let
\begin{equation}
T_{\alpha}=\{x:\:  |P_{\alpha}(x)|<d_{\alpha},\,\langle x,v_{\alpha}\rangle>R_{\alpha}\},
\end{equation}
where
\begin{equation}
P_{\alpha}(x)=x-\langle x,v_{\alpha}\rangle v_{\alpha}
\end{equation}
is the orthogonal projection onto the hyperplane $\{w: \langle w, v_\alpha \rangle = 0\}$.
We assume here that $\epsilon\in\cC^{\infty}(\bbR^d),$ and for each
$\alpha\in\cA$,
$\epsilon_{\alpha}(x)\overset{d}{=}\epsilon(x)\restrictedto_{T_{\alpha}}.$ It is also
assumed that $\epsilon_{\alpha}(x)$ depends only on $P_{\alpha}(x)$ for $\langle
x,v_{\alpha}\rangle>R'_{\alpha} > R_\alpha$, see Figure~\ref{WaveGuides}.

To rephrase this as a scattering problem we let $k_1^2=\epsilon_1f^2,$ define
the potential $q(x)=(\epsilon(x)-\epsilon_1)f^2,$ and replace~\eqref{eqn3.200}
with
\begin{equation}\label{eqn4.200}
(\Delta+q(x)+k_1^2)u=0.
\end{equation}
The potential $q(x)$ is supported in the set $\Omega$ defined above.

In time--independent scattering theory, one imagines that there is an `incoming'
field $u^{\In}(x),$ which is constructed from either a wave-guide mode (see
Sections~\ref{sec3} and~\ref{sec.lap_scat}) for one of the channels, or a free-space wave packet (see
Section 6 of~\cite{EpWG2023_1}).  We then look for an `outgoing' solution to
\begin{equation}\label{eqn5.220}
  (\Delta+q(x)+k_1^2)u^{\out}=-(\Delta+q(x)+k_1^2)u^{\In},
\end{equation}
so that $u^{\tot}=u^{\In}+u^{\out}$ represents the total field that results from
the incoming field scattering off of the wave-guide network.

There are many papers in the Applied Math and Physics literature that consider
such problems;
see~\cite{BonnetBendhia_etal,BonnetTillequin2001,ChandlerMonkThomas2007,
  ChandlerZhang} in the Applied Math literature,
and~\cite{Nosich1994,KNH_SIAM_2005} and the references therein for the
Mathematical Physics literature.  What these various papers do not provide,
however, are rigorous, physically motivated definitions of the concepts
`incoming' and `outgoing' for the full $d$-dimensional problem.  \Rd
In~\cite{CiraoloMagnanini2008} the problem of $2d$ bi-infinite wave-guides is
treated and an outgoing radiation condition is given that implies uniqueness of
the solution. While their analysis and radiation condition are limited to this
special case, solutions satisfying their condition also satisfy the
condition we present herein. \Bk In general, earlier works do not prove a uniqueness
result for this scattering problem, which is a central component of a complete
theory. The fact that the potential $q(x)$ does not vanish at infinity is what
makes this challenging, and not covered by the more standard `two-body'
scattering literature.  In particular, the classical Sommerfeld radiation
conditions must be reformulated to give new and tractable criteria for
uniqueness in this setting.

\begin{figure}
  \centering \includegraphics[width= 10cm]{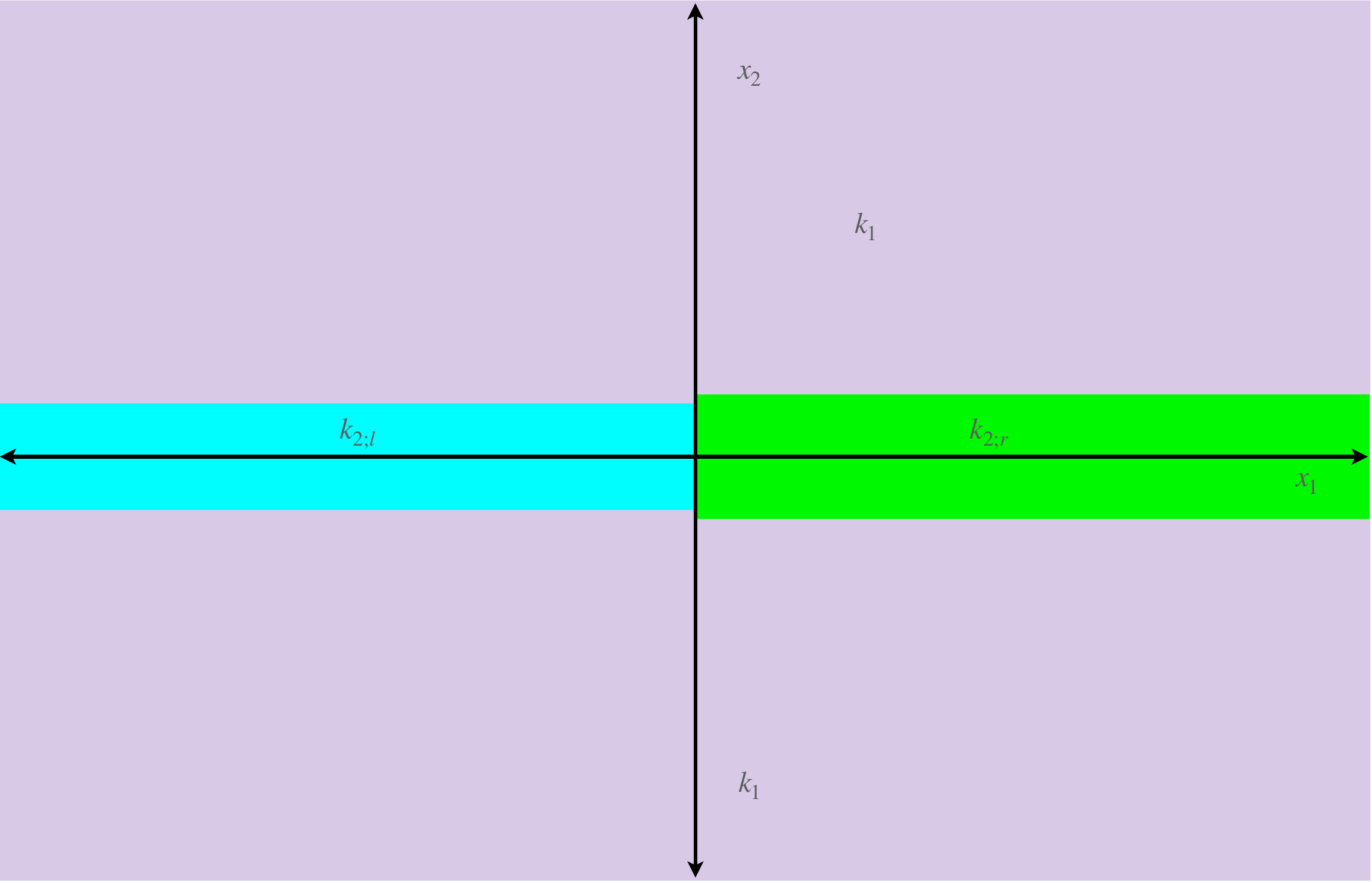}
    \caption{Two dielectric channels in $\bbR^2$ meeting along a straight interface.}  
   \label{fig0}
\end{figure}

\begin{figure}
  \centering \includegraphics[width= 14cm]{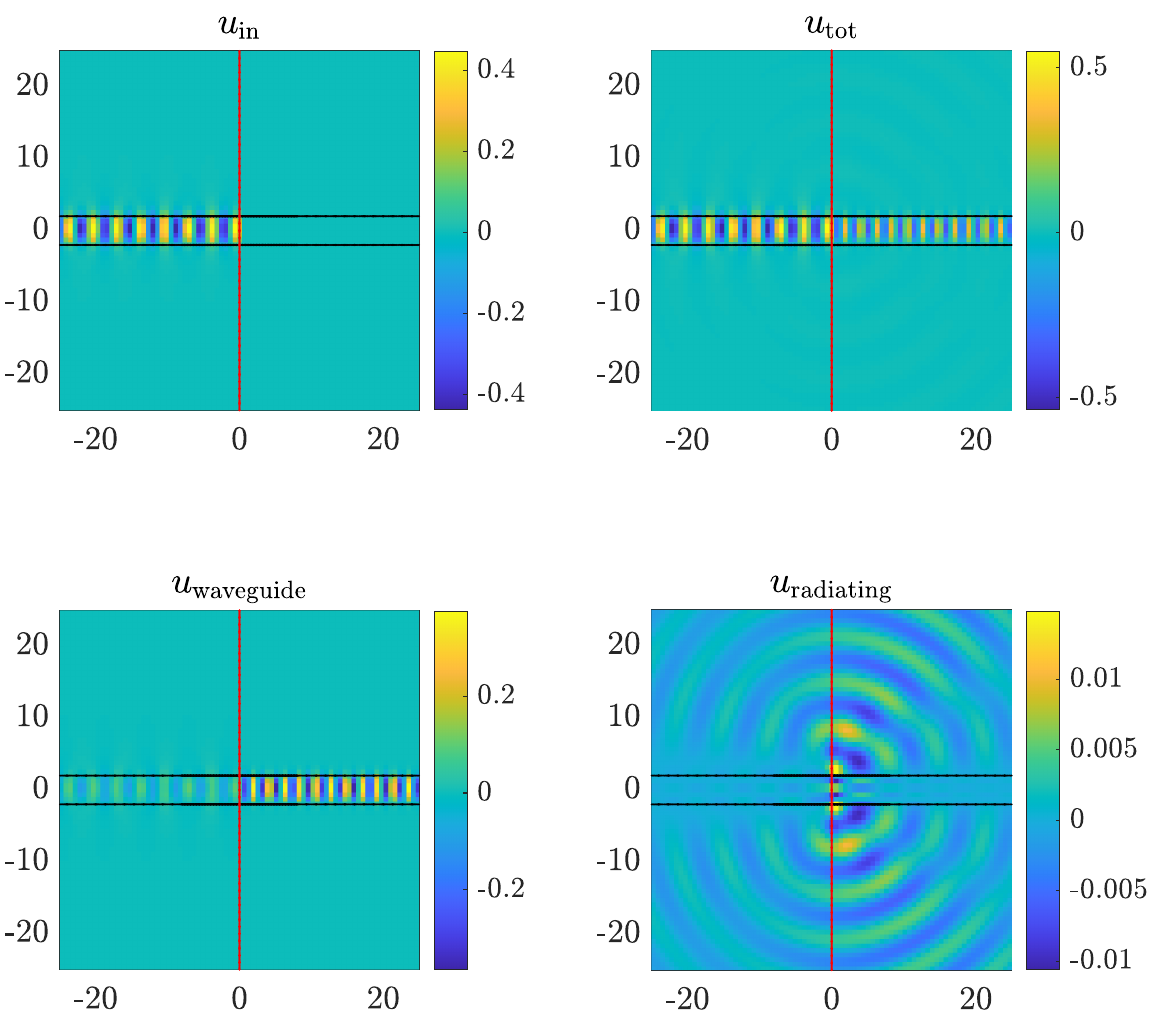}
    \caption{The upper left panel shows the incoming waveguide mode, $u^{\In},$
      the upper right panel shows the total field, $u^{\tot}.$ The lower left
      panel shows the wave-guide mode part of $u^{\out},$ and the lower right panel
      the radiation part of $u^{\out}.$ These numerical solutions and plots are
      obtained using MATLAB software written by Tristan Goodwill and described
      in~\cite{GE_2024}.}
   \label{scat_field}
\end{figure}

To the best of our knowledge, none of the approaches in the literature for
solving this problem are amenable to a numerical realization that provides an
accurate representation of the radiation field outside of the channels.  This
latter issue is addressed in~\cite{EpWG2023_1,EpWG2023_2} for the simple model
problem of two semi-infinite rectangular channels in $\bbR^2$ that meet along a
common perpendicular line, with a piecewise constant potential, see
Figure~\ref{fig0}.  These papers reformulate the scattering problem as a
transmission problem, which is solved using an integral equation approach. This
representation also yields precise asymptotics for the solutions.  A numerical
method of solution using this formulation is given
in~\cite{GE_2024}. Figure~\ref{scat_field} shows a numerical example of the
solution to such a scattering problem, using a wave-guide network model like that in
Figure~\ref{fig0}.

In the present paper we provide radiation conditions that imply
uniqueness.  We assume here, to conform with earlier literature, that the
potential, $q(x),$ is smooth.  The integral equation method
from~\cite{EpWG2023_1,EpWG2023_2} can easily be adapted to handle smooth
potentials, see~\cite{GHRQ}. On the other hand, by a modification of the
techniques used here, such as those in~\cite{PSS1981}, it should be possible to
extend the results of this paper to the piecewise constant case.

In fact, radiation conditions for similar multi-channel scattering problems were
treated in the mathematical physics literature as far back as the
mid-1990's. In~\cite{Isozaki94} H.~Isozaki provided a solution to this problem
in the setting of $N$-body Schr\"odinger potentials. The assumptions he
makes on the potential exclude the ones considered here. Nonetheless his
results apply to open wave-guide networks as shown by a different analysis of the
problem due to Vasy, culminating in the paper~\cite{VasyAsterisque}.  Vasy uses
the tools of geometric microlocal analysis, significantly refining a
method introduced by Melrose in~\cite{Melrose94} to study scattering theory for
potentials with sufficient decay.  The key new idea proposed by Melrose is that
the oscillations `at infinity' and decay of solutions to $(\Delta + q(x) + k_1^2)u = 0$ can be
encoded in a geometric notion called the {\it scattering wave-front set}. He
shows that this set is strongly constrained by the asymptotic behavior of the
operator itself.  Vasy's work generalizes this approach considerably,
incorporating potentials that do not decay along certain (possibly higher
dimensional) channels.

The principal goal of the present paper is to present the Isozaki/Vasy radiation
conditions in a form usable for open wave-guide networks.  Since 
the channels are concentrated along rays (rather than higher dimensional
subspaces), many details in~\cite{Isozaki94,Melrose94,VasyAsterisque} simplify.
The work of Melrose and Vasy also gives a precise explanation for the nature of the radiation
conditions, and explains why Isozaki's conditions suffice.  We give a brief
accounting of their work, but refer to those papers for the detailed
technical arguments. We also apply these conditions to the solutions obtained
in~\cite{EpWG2023_1,EpWG2023_2}, showing that they are indeed outgoing and
agree with the so-called `limiting absorption solution,' whose existence is demonstrated
in~\cite{VasyAsterisque}.

One important aspect consequence of this analysis is that radiation conditions
for open wave-guide networks are essentially local at infinity.  The following
language makes this precise. Denote by $\overline{\bbR^d}$ the radial
compactification of $\bbR^d,$ which adds a boundary point for each direction.
Altogether we add a unit sphere $S^{d-1}$ consisting of all such asymptotic
directions, and thus $\overline{\bbR^d}$ is diffeomorphic to the closed unit
ball.
It is convenient to start with polar coordinates, $r = |x| \geq 0$, $\omega = x/|x| \in
S^{d-1}.$ 
Setting $\rho = 1/r$ gives coordinates, $(\rho,\omega),$  near
$\pa\overline{\bbR^d}$ with
\begin{equation}
\pa\overline{\bbR^d}=\{\rho=0,\,\omega\in S^{d-1}\}.
\end{equation}

The points $\cC=\{v_{\alpha}:\:\alpha\in\cA\}\subset\pa\overline{\bbR^d},$
corresponding to the directions of the channels, are called the \emph{channel
ends}, and the complement $\cC^c = S^{d-1}\setminus \cC$ is called the
\emph{free boundary}.  We localize a function near a point $\omega_0 \in \pa
\overline{\bbR^d}$ by multiplying by a cut-off function $\Psi \in \mathcal
C^\infty(\overline{\bbR^d})$ which equals $1$ near $\omega_0$ and vanishes
outside a small neighborhood of that point.  For example, we can use smooth
`conical cut-offs'
\begin{equation}
\Psi_{\omega_0,\delta}(x)=
  \begin{cases}
    &1\text{ for }\left|\frac{x}{|x|}-\omega_0\right|<\delta, |x|>C+1,\\
    &\phantom{=}\\
    &0\text{ for }\left|\frac{x}{|x|}-\omega_0\right|>2\delta, |x|<C,    
  \end{cases}
\end{equation}
for any $C > 0$; the parameter $\delta >0$ measures the aperture of the conical neighborhood.   In the sequel, the operation
\begin{equation}
u \longmapsto \mathcal F( \Psi_{\omega_0, \delta} u)
\end{equation}
is sometimes referred to as a {\it conic Fourier transform}.

Isozaki defines the notion of incoming ($+$) and outgoing ($-$) solutions by using  certain classes of \emph{admissible} 
operators $\cR^{\ka}_{\pm}(\epsilon)$: 
a solution is declared to be incoming, resp. outgoing, near $\omega_0$ if there there exists $\ka_0,\epsilon>0$ and $\delta,\delta'> 0$  such that for any $P_{\pm}\in\cR^{\ka_0}_{\pm}(\epsilon)$
\begin{equation}
\rho^{\frac 12-\eta}\Psi_{\omega_0,\delta'}P_{\pm}[\Psi_{\omega_0,\delta}u]\in L^2(\bbR^d),
\end{equation}
for some $\eta > 0.$ Notice that $\rho^{\frac12 - \eta} f \in L^2$ means that
(at least on average) $f$ decays faster than $r^{\frac{1-d}{2}}$.  The operators
in $\cR^{\ka}_{\pm}(\epsilon)$ are essentially pseudodifferential operators,
which are described in Remark~\ref{rmk7.205}.  Unlike the classical Sommerfeld radiation operators,
$\pa_r\mp ik_1,$ they do not depend, in any essential way, on the background
wave number, $k_1^2.$ The Melrose/Vasy theory explains why this is not
necessary. 

It is clear from Isozaki's treatment that the radiation conditions
reflect the singularities of the conic Fourier transforms of the
solution $u.$ This connection is made explicit in the classical
setting of potentials vanishing at infinity in~\cite{Melrose94} and
leads directly to what Melrose calls the \emph{scattering wave-front}
set, as mentioned above. Vasy~\cite{VasyAsterisque} extends this
notion to a class of non-compactly supported potentials that includes
open wave-guide network potentials. Isozaki's admissible operators are
subsumed in a larger and more flexible class of pseudodifferential
operators, which constitute the \emph{3-body scattering calculus}.
The ensuing microlocal treatment explains many aspects of the
radiation conditions in terms of propagation phenomena, at infinity,
analogous to those that arise in the analysis of hyperbolic
equations. Vasy also proves that a solution to
$(\Delta+q+k_1^2)u=f\in\cS(\bbR^d),$ which is either incoming or
outgoing, in this extended sense, is unique.  This is used to show the
existence of the limiting absorption resolvents $(\Delta+q+k_1^2\pm
i0)^{-1}=\lim_{\sigma\to 0^+}(\Delta+q+k_1^2\pm i\sigma)^{-1}$, as
bounded maps between weighted Sobolev spaces, which produce the
unique incoming ($-$) and outgoing ($+$) solutions.

We remark here that both Isozaki and Vasy work in a setting that includes more
general $N$-body potentials in Euclidean space. One key new feature, when $N >
3$, is that the channels (which lie along positive dimensional subspaces and not just
lines) may intersect at infinity; this does not occur when $N=3$.  Our geometric
assumptions about channels ensure that this does not happen, and for this reason
we may employ tools adapted to the setting of 3-body potentials.  Moreover, in most of the
$N$-body Schr\"odinger equation literature the channels are at least 2
dimensional and so they meet $\pa\overline{\bbR^d}$ in a positive dimensional
set. This precludes the use of conic localization near the ends of the channels,
which is a significant simplification in the wave-guide network  case.

We begin in Section~\ref{sec2} by recalling the classically understood behavior
of solutions to $(\Delta+q(x)+k_1^2)u\in\cS$ where the potential $q(x)$ decays at
infinity. In Section~\ref{sec3} we give a more detailed description of the
wave-guide network models we consider, as well as the definition of wave-guide
modes. In Section~\ref{sec.rad_cond} we describe two versions of the
radiation conditions for  open wave-guide networks. Isozaki's condition
requires less background to explain, but to present Vasy's  conditions,
we must first describe both Melrose's scattering calculus and Vasy's 3-body scattering
calculus.  In Section~\ref{sec_2d} we prove that the
solutions obtained in~\cite{EpWG2023_1,EpWG2023_2} are outgoing in the sense of
Isozaki/Vasy. These calculations illustrate our assertion that Vasy's
conditions are more flexible and easier to check in specific examples.  We
finally prove the missing uniqueness theorem for the Fredholm equations
of index zero introduced in~\cite{EpWG2023_1}, which thereby completes the proof of the
existence of solutions to the transmission problem, see~\eqref{eqn99.221}, \eqref{eqn88.220}.

Additional results are given in two appendices. In the first we show, by direct
computation, that a classically outgoing solution satisfies Isozaki's form of
the radiation condition. In the second we show that the channel-to-channel
scattering coefficients are well defined for an open wave-guide network.
\smallskip

\centerline{\bf Notational Conventions}
\begin{enumerate}
\item To simplify notation the letter $k,$ without an argument, refers to the
  background wave-number, $k^2=\epsilon_1 f^2.$
  \item When it will not cause confusion, we refer to the image of the
    projections $\{P_{\alpha}:\:\alpha\in\cA\}$ as $\bbR^{d-1}.$
  \item In Section~\ref{sec_2d}, and Appendix~\ref{App1}, $l$ corresponds to $x_1<0,$ and $r$ to $x_1>0.$
     \item In Section~\ref{sec_2d}, and Appendix~\ref{App1}, $g_k(x-y)$ is the
       kernel for the outgoing fundamental solution,   $(\Delta+k^2+i0)^{-1},$ in 2-dimensions.
       \item In Section~\ref{sec_2d}, and Appendix~\ref{App1}, $\cS_k$ and
         $\cD_k$ are the single and double layer potentials over the $x_2$-axis
         defined by $g_k.$
\end{enumerate}

\smallskip

{\small
\centerline{\bf Acknowledgments}

\noindent
The authors wish to thank Andras Vasy for very helpful conversations at various
stages in the preparation of this project, and for his careful reading and
comments on an earlier version of this paper. We would also like to thank the
referee for his/her careful reading and many  suggestions that substantially
improved our paper.}

\section{Channel-free scattering}\label{sec2}
Before turning to the main topic of interest here, we briefly review some
classical facts about scattering in Euclidean backgrounds.  We consider both the
`free' case, i.e. for the Helmholtz equation $(\Delta + k^2) u = 0$ in $\mathbb
R^d$, and then a few generalizations to the equation $(\Delta + q(x) + k^2) u =
0$ where $q(x)$ is a smooth, short range potential, i.e. vanishing at infinity
at least as fast as $|x|^{-1-\epsilon}$ for some $\epsilon > 0$.  All of this
material is classical and well-known.

Consider first the space of solutions to $(\Delta + k^2) u = 0$ in all of
Euclidean space.  There are two key building blocks for all other solutions: the
first is the class of plane wave solutions:
\[
u_{k,\omega_0} (x) = e^{ik x \cdot \omega_0},
\]
where $\omega_0 \in S^{d-1}$ is fixed. While these seem quite simple, they are
actually singular in that they they do not decay at infinity.  There is another
basic set of solutions; in $\bbR^2$ these are the functions
\[
 Z_\ell (kr) e^{\pm i \ell \phi}, \quad k, \ell \in \mathbb N_0,
 \]
 with $Z_{\ell}$ any Bessel function of order $\ell.$

 These have analogues in dimensions $d>2,$ where the factors $e^{\pm i \ell
   \phi}$ are replaced by spherical harmonics: For any solution to
 $\Delta_{S^{d-1}}Y=-\lambda^2Y,$ there are solutions to the Helmholtz equation
 of the form $r^{\frac{2-d}{2}}Z_{\nu}(kr)Y(\omega)$.  Here $Z_{\nu}(r)$ is any
 Bessel function of order
 $\nu={\pm\sqrt{\lambda^2+\left(\frac{d-2}{2}\right)^2}}.$ As the eigenvalues
 on $S^{d-1}$ are $\lambda^2 = \ell (d-2 + \ell)$, $\ell \in \mathbb N_0$, the
 degrees $\nu_\ell = \pm\left( \ell+\frac{d-2}{2}\right)$.  Using standard
 asymptotics of Bessel functions, we have that the globally smooth solutions in
 $2d$ satisfy
$$ v_{k,\ell}^{\pm}(r,\phi)\overset{d}{=}J_\ell (kr) e^{\pm i \ell \phi} \sim \frac{ e^{\pm
     i\ell \phi}}{\sqrt{r}}( a_0^+e^{ikr} + a_0^- e^{-ikr}),\text{ as }r\to\infty,
$$
for some constants $a_0^\pm.$

There appears to be a stark difference between these two classes of solutions:
the $u_{k,\omega_0}$ are bounded but do not decay, and they oscillate, but with `radial
frequency' reaching a maximum only in the directions $\pm \omega_0$.  By
contrast, the $v_{k,\ell}^\pm$ decay at a fixed rate, and oscillate
uniformly in the radial direction at frequencies $\pm k.$ These solutions are closely related, by virtue
of the classical formula
\[
J_\ell(kr) = \int_{S^{1}} \left[A_\ell^+(\omega) e^{ik\omega \cdot x} + A_\ell^-(\omega) e^{-ik\omega \cdot x}\right]\, dV_\omega,
\]
for certain smooth functions $A_\ell^\pm(\omega)$ on the circle.  This formula
is basically equivalent to the classical integral definition of the $J$-Bessel
function:
\[
J_\ell(r) = \frac{i^{-\ell}}{\pi}\int_0^\pi e^{ir \cos \phi} \cos (\ell\phi)\, d\phi.
\]
These smooth, decaying solutions are simply averages over all directions of the
plane wave solutions, with suitably chosen amplitudes.

Motivated by these two examples, we quote a result about general solutions to $(\Delta + k^2) u = 0$ in any external
region $|x| \geq R$. Namely, any such solution admits an asymptotic expansion of the form
\begin{equation}
u(x) \sim  r^{\frac{1-d}{2}} \left[e^{ikr} \sum_{j=0}^\infty a_j^+(\omega)r^{-j} + e^{-ikr} \sum_{j=0}^\infty a_j^-(\omega)r^{-j}\right].
\label{expansion}
\end{equation}
The meaning of this asymptotic expansion is unexpectedly subtle in that the
coefficients $a_j^\pm(\omega)$ are, in general, only distributions on the
sphere.  Thus to make sense of this, we must average, i.e., integrate against a
smooth test function in $\omega$.  Thus the actual meaning of \eqref{expansion}
is that if $\chi(\omega) \in \mathcal C^\infty(S^{d-1})$ is arbitrary, then
\begin{equation}
\begin{split}
& \int_{S^{d-1}} u(r\omega) \chi(\omega)\, dV_\omega  \\ 
& \qquad \qquad \sim r^{\frac{1-d}{2}}\left[  e^{ikr}  \sum_{j=0}^\infty\langle a_j^+, \chi \rangle r^{-j} 
+ e^{-ikr} \sum_{j=0}^\infty \langle a_j^-, \chi \rangle r^{-j}\right].
\end{split}
\label{expansion2}
\end{equation}
\Rd Here $\langle\cdot,\cdot\rangle$ is the canonical pairing
$\cC^{-\infty}(S^{d-1})\times \cC^{\infty}(S^{d-1})\to\bbC.$ \Bk
As a function of the radial variable $r$ alone, this is an asymptotic expansion,
in the traditional sense, for a function of one variable.  The existence of such
expansions, when the leading coefficients are smooth functions on $S^{d-1},$ is
sketched in \cite[Section 1.3]{MelroseGST} using stationary phase. The existence
of a leading distributional term in general is also noted there; the complete
distributional expansion can be obtained by an iterative argument using the
Mellin transform in the radial variable.

The functions $v_{k,\ell}$ satisfy these asymptotic conditions in the usual
strong sense since the amplitude functions $a^\pm_0e^{\pm il\phi}$ are smooth. The
plane wave solutions admit an expansion of this form, albeit with distributional
coefficients. Writing $x = r\omega$, we have
\begin{equation}\label{eqn14.210}
e^{ik r\omega \cdot \omega_0} \sim  c_dr^{\frac{1-d}{2}} [e^{ikr-i\pi/4}
  \delta_{\omega_0}(\omega) + e^{-ikr+i\pi/4} \delta_{-\omega_0}(\omega)] +O(r^{-\frac{d+1}{2}}).
\end{equation}
To see this, we apply stationary phase, which shows that, as $r\to\infty,$
\begin{multline}
\int_{S^{d-1}} e^{ikr \omega \cdot \omega_0} \chi(\omega) \, dV_\omega =\\  c_dr^{\frac{1-d}{2}}[e^{ikr-i\pi/4} \chi(\omega_0) + e^{-ikr+i\pi/4} 
\chi(-\omega_0)]+O(r^{-\frac{d+1}{2}}),
\end{multline}
for any $N>0.$



\begin{definition}\label{def1}
A formal solution $u\in\cS'(\bbR^d)$ to $(\Delta + k^2) u = f\in
\mathcal S(\mathbb R^d)$ is called outgoing if, in \eqref{expansion},
all of the coefficients $\{a_j^-(\omega):\: 0\leq j\}$ vanish
identically. In other words, a solution is said to be outgoing if
\[
u(r \omega) \sim r^{\frac{1-d}{2}}e^{ikr} \sum_{j=0}^\infty r^{-j} a_j^+(\omega).
\]
\end{definition}
It is a non-trivial fact, proved in~\cite{Melrose94},  that if $u$ is outgoing (and $f$ is
Schwartz), then all of the coefficients $a_j^+(\omega)$ are $\mathcal C^\infty$.
Granting this, then an integration by parts leads to the conclusion
that if $f \equiv 0$ and $u$ is outgoing, then $u \equiv 0$. In other words,
there are no non-trivial, globally defined outgoing solutions, see~\cite{Rellich1943}.  There is an
analogous notion that a solution is called {\it incoming} if all of the `positive'
coefficients, $\{a_j^+(\omega)\},$ in its expansion vanish.  Everything that we say here about
outgoing solutions has a direct analogue for incoming solutions.

The explicit solutions $v_{k,\ell}$ we wrote down earlier are neither incoming
nor outgoing, since both coefficients $a_0^\pm$ are non-vanishing.  However,
there is another solution of the Bessel equation called the Hankel function of
the first kind, denoted $H_\nu^{(1)}$, so that $w_{k,d} (r\omega)=
r^{\frac{2-d}{2}} H_{(d-2)/2}^{(1)}(kr) $ is an outgoing solution to the
Helmholtz equation in $\bbR^d\setminus \{0\}.$ It has a singularity at $r = 0$, which is
in agreement with the fact that there are no globally defined outgoing
solutions. This solution is used to define the outgoing fundamental solution
for $(\Delta+k^2).$

The classical Sommerfeld condition provides a criterion to check whether a
solution $u$ is outgoing without the need for an asymptotic expansion. Observe
that even in the `best' case where all $a_j^\pm$ are smooth, the leading terms
$e^{\pm ikr} r^{(1-d)/2}$ do not lie in $L^2(\bbR^d)$.  However, applying the operator
$\partial_r - ik$ to $u$ reduces the order of growth of one of these terms:
\[
\left(\partial_r - ik \right) e^{ikr} r^{\frac{1-d}{2}} =  \left(\frac{1-d}{2} \right)e^{ikr} r^{-\frac{1+d}{2}},
\]
which {\it does} lie in $L^2(\bbR^d)$.  On the other hand, applying $\partial_r - ik$ to $e^{-ikr} r^{(1-d)/2}$ does not yield
better decay. As the coefficients $a_j^{\pm}$, $j > 0$, all depend linearly on
$a_0^{\pm}$, it follows that
\begin{equation}
\mbox{The solution}\ u\ \mbox{is outgoing if}\ \left(\partial_r - ik \right) u \in L^2(\mathbb R^d).
\label{SF}
\end{equation}
In fact the weaker conditions
\begin{equation}\label{eqn13.221}
  r^{\delta-\frac 12}(\pa_r-ik)u(r\omega)\in L^2(\bbR^d)\text{ for a }\delta>0
\end{equation}
suffice to conclude that a solution is outgoing.

Note that given any $f \in \mathcal S$, there exists a unique
outgoing solution $u$ to $(\Delta + k^2) u = f$.  This solution is obtained by
convolving with the outgoing resolvent:
\[
u = \int_{\mathbb R^d} R_{k^2+i0}(x , \tilde{x}) f(\tilde{x}) \, dV_{\tilde{x}}.
\]
The integral kernel of the outgoing resolvent is the radial function
\begin{equation}
  R_{k^2+i0}(x , \tilde{x})=C_d\frac{H^{(1)}_{\frac{d-2}{2}}(k|x-\tx|)}{|x-\tx|^{\frac{d-2}{2}}}.
\end{equation}

Everything we have said here has an analogue for   general Schr\"odinger operators
\[
(H+k^2)u := (\Delta + q(x) + k^2) u,
\]
where $q $ is sufficiently regular and decays sufficiently rapidly at
infinity. A solution to $(H+k^2)u=f\in\cS(\bbR^d)$ is outgoing if it
satisfies~\eqref{eqn13.221}. We no longer have explicit solutions, but one can
still prove the existence of an outgoing resolvent $R(k^2+i0)$, as well as the
existence of expansions for solutions to $(H+k^2)u = f$.  The
existence of the outgoing resolvent is called the \emph{limiting absorption
principle}: If $\sigma>0,$ then the resolvent operators $R(k^2\pm
i\sigma)=(H+k^2\pm i\sigma)^{-1}$ are well defined as bounded
operators on $L^2(\bbR^d),$ which, in fact map $\cS(\bbR^d)$ to itself. The
limiting absorption principle states that,
\begin{equation}
R(k^2\pm i 0)=  \lim_{\sigma\to 0^+}R(k^2\pm i\sigma)
\end{equation}
exist as bounded maps from $r^{-(\delta+\frac 12)}L^2(\bbR^d)$ to
$r^{\delta+\frac 12}L^2(\bbR^d),$  for any $\delta>0.$ Moreover, if
$f\in\cS(\bbR^d),$ then
\begin{equation}
  u_{\pm}=R(k^2\pm i 0)f
\end{equation}
is the unique outgoing ($+$), resp. incoming ($-$) solution to
\begin{equation}
  (\Delta + q + k^2) u_{\pm}=f.
\end{equation}
These extensions are covered in~\cite{Melrose94} and~\cite{TaylorPDEII}.

\section{Open Wave-Guide Networks and Wave-Guide Modes}\label{sec3}
\begin{figure}
  \centering
  \includegraphics[width= 8cm]{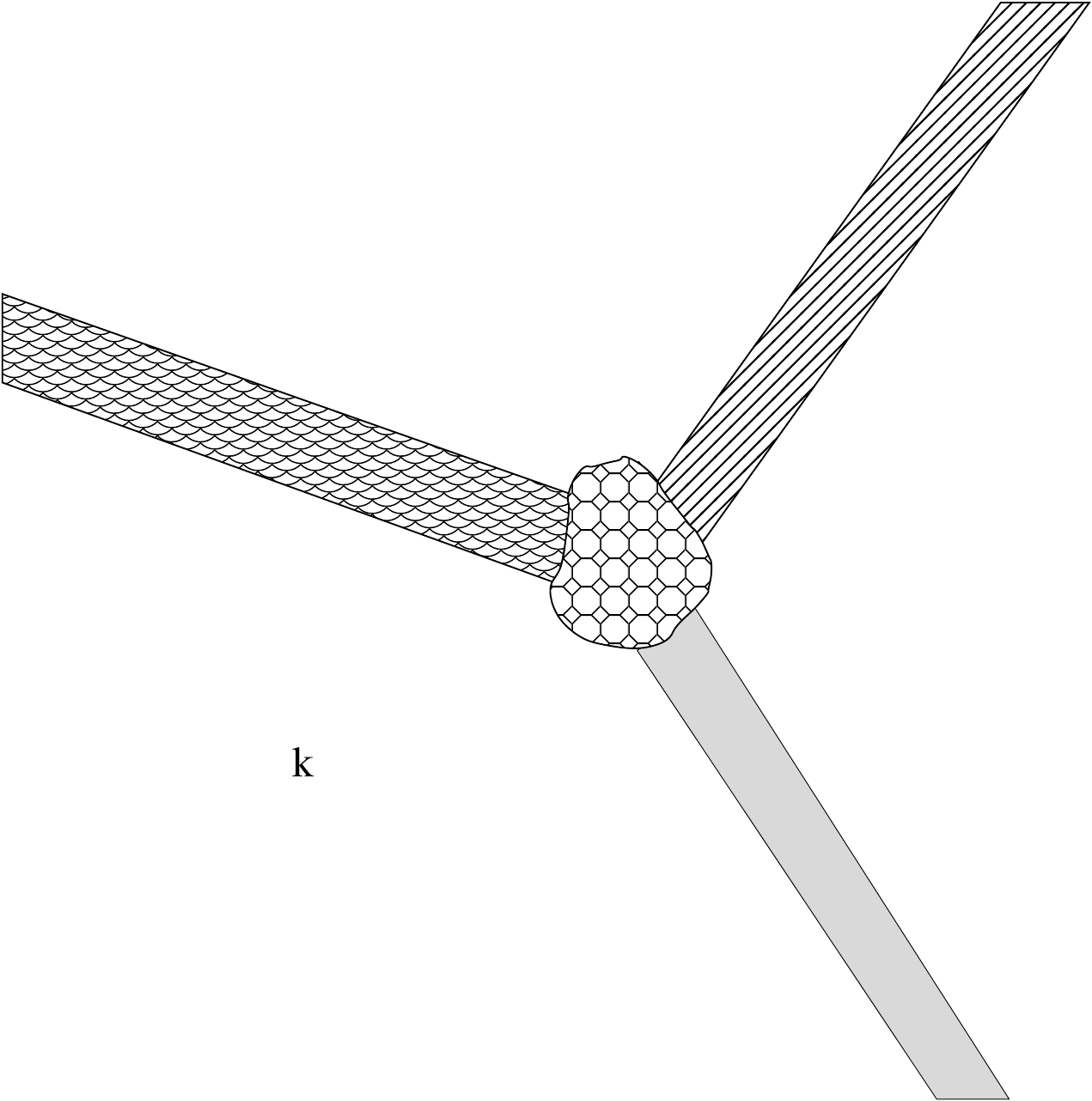}
    \caption{Three dielectric channels, indicated as tubes around rays
      extending to infinite, meeting in a compact interaction zone, $\Omega_0.$
    }
    \label{WaveGuides}
\end{figure}

The main goal of this paper is to provide physically motivated radiation
conditions that imply uniqueness for open wave-guide networks. In this section we
describe the mathematical model we use for an open wave-guide network  in $\bbR^d:$ as
described in Section~\ref{sec1}, imagine a collection of dielectric,
non-conducting, non-magnetic, `pipes,' which interact in a compact region of
space, and are asymptotic to disjoint straight lines as they head off to
infinity, see Figure~\ref{WaveGuides}. Within the pipes, contained in $\Omega,$  the permittivity is
spatially dependent, whereas in the exterior region it assumes a constant
positive real value, $\epsilon_1.$  As the solution to the underlying wave
equation,  \eqref{eqn3.201}, is time harmonic, taking the
form $e^{-ift}u(x),$ the `free space' wave number is $k^2=\epsilon_1f^2.$

We next describe the channel regions. As in the introduction, fix a finite set of asymptotic directions
$\cC=\{v_{\alpha}:\:\alpha\in\cA\}\subset\pa\overline{\bbR^d}\simeq S^{d-1}$, and for each $\alpha$,
denote by $P_{\alpha}$ the orthogonal projection onto the orthocomplement of $v_\alpha$.
Following standard conventions in $N$-body scattering theory, we write
\begin{equation}\label{eqn18.222}
x_\alpha = \langle x, v_\alpha\rangle,\ \ \mbox{and}\ \ \ x^\alpha = P_\alpha x = x - x_\alpha v_\alpha.
\end{equation}

For $\alpha\in\cA,$ let $\epsilon_{\alpha}\in\cC^{\infty}(\Im P_{\alpha})$ be
a positive, real valued
function that describes  the permittivity within the $\alpha^{\Th}$ channel for
$\langle x, v_\alpha\rangle$ sufficiently large. At frequency $f$ we set
$$q_{\alpha}(x)=(\epsilon_{\alpha}(x^{\alpha})-\epsilon_1)f^2;$$
by assumption
$q_{\alpha}$ is smooth and has compact support as a function of $x^{\alpha}.$ This represents
the deviation from the free space wave-number within this channel. Finally, choose
$\psi_+\in\cC^{\infty}(\bbR),$ vanishing for $t < 1,$ and equals to $1$ for $t >
2$. There is a compactly supported function $q_0\in\cC^{\infty}_c(\bbR^d)$ and positive constants
$\{r_{\alpha}\},$  so that
\begin{equation}
q(x)=(\epsilon(x)-\epsilon_1)f^2=q_0(x)+\sum_{\alpha\in\cA}q_{\alpha}(x^\alpha)\psi_+(r_{\alpha} x_\alpha ).
\label{formq}
\end{equation}

In terms of this data, the model for the open wave-guide network at
frequency $f$ is given by the time-independent Schr\"odinger operator
\begin{equation}\label{eqn93.211}
H+k^2=(\Delta+q+k^2).
\end{equation}
The support of $q(x)$ outside a 
compact set lies in a bounded neighborhood of the collection of rays $\{r
v_{\alpha}:\: r\gg 0,\,\alpha\in\cA\}.$ These meet $\pa\overline{\bbR^d} = S^{d-1}$
at the points
\[
\cC=\{v_{\alpha}:\: \alpha\in\cA\}.
\]
We call these points the `channel ends.' From the perspective of $N$-body Schr\"odinger
operators, this is called the 3-body case because the channel ends are disjoint.

For simplicity we have assumed that the wave number within a wave-guide channel is
independent of $x_\alpha$ when $x_\alpha \gg 0$. It is not difficult to
adapt the methods  below to allow $q\restrictedto_{T_\alpha}$ to decay, at some
sufficiently high rate, to an $x_\alpha$-independent function of $x^\alpha$ as
$x_\alpha \to \infty$. 

The range of $P_{\alpha}$ is a hyperplane in $\bbR^d;$ denote by
$\Delta_{\alpha}$  the Laplace operator in an orthonormal coordinate system on
the range of $P_\alpha.$ For any $\alpha$,  we can express
\begin{equation}
\Delta=\pa_{x_{\alpha}}^2+\Delta_{\alpha}.
\end{equation}
 Hence in a conic neighborhood of the $\alpha^{\mathrm{th}}$ channel,  where
$x_{\alpha}$ is sufficiently large, 
\begin{equation}
H = \pa_{x_\alpha}^2 + H_\alpha, \quad \mbox{where}\qquad H_\alpha = \Delta_\alpha + q_\alpha(x^\alpha).
\end{equation}
The $\{H_\alpha\}$ are called the subsystem Hamiltonians. 

Our interest is in \emph{formal solutions}, i.e., tempered distributions $u$
satisfying $(H+k^2)u\in\cS(\bbR^d)$.  By classical ellipticity, any formal
solution is in $\cC^{\infty}(\bbR^d)$ and of polynomial growth at infinity.  As
a first step toward describing the possible asymptotics of formal solutions as
$|x| \to \infty$, we introduce special solutions associated to each
subsystem. If $q_{\alpha} \not\equiv 0,$ then there may be a finite set of
$E\in\bbR\setminus\{0\},$ such that the space
of $L^2$-solutions to
\begin{equation}
H_\alpha u(x^\alpha) = E^2u(x^{\alpha})
\end{equation}
is non-trivial. Any such solution necessarily belongs to $\cS(\bbR^{d-1}).$  Note
that we have excluded $E=0$ in order to avoid various technical issues. 

At the risk of proliferation of
subscripts, let $0<E_{\alpha,1}^2 < \ldots < E_{\alpha, N_\alpha'}^2$ be the list of values for which this eigenspace
is nontrivial, and for each $E_{\alpha, \ell}$, let $\{u_{\alpha, j, \ell}:\ j = 1, \ldots, N_{\alpha,\ell}''\}$ be a fixed orthonormal
basis for this eigenspace. We then have that 
\begin{equation}\label{eqn23.200}
(\pa_{x_{\alpha}}^2+\Delta_{x^{\alpha}}+q_{\alpha}(x^{\alpha})+k^2)[e^{\pm  ix_{\alpha}\sqrt{E^2_{\alpha,\ell}+k^2}}u_{\alpha,j,\ell}(x^{\alpha})]=0.
\end{equation}
The functions within the brackets  in~\eqref{eqn23.200} are called \emph{wave-guide modes}; they do
not decay as $x_{\alpha}\to\infty,$ but are strongly localized within the
$\alpha^{\Th}$ channel.  With respect to the time dependent factor $e^{-if
  t},\,f>0,$ it is then reasonable (and consistent) to call the solution with
exponential factor $e^{ix_{\alpha}\sqrt{E^2_{\alpha,\ell}+k^2}}$ outgoing, and
the solution with $e^{-ix_{\alpha}\sqrt{E^2_{\alpha,\ell}+k^2}}$ incoming.

We can construct global formal solutions out of these wave-guide modes. Indeed,
suppose that $u_\alpha$ is a wave-guide mode for the channel $\cC_\alpha,$ with
energy $E_{\alpha,\ell}.$ Choose
a smooth cutoff function, $\chi_{\alpha},$ which equals $1$ in a small conical
neighborhood centered on the ray
$\{rv_\alpha:\: r\gg R\},$  and  vanishes outside a slightly larger conical neighborhood
around this axis. The larger cone is chosen so that $\supp\chi_{\alpha}$ is disjoint from
$\supp[q_{\beta}(x^{\beta})\psi_+(r_{\beta}x_{\beta})],$
  for $\beta\neq\alpha.$   With these choices we let
\[
(H +k^2)w= f, \text{ where }w= \chi_\alpha(x_\alpha,
x^\alpha) e^{\pm ix_{\alpha}\sqrt{E^2_{\alpha,\ell}+k^2}}u_\alpha(x^\alpha).
\]
From the exponential decay of  $u_\alpha$ as $|x^\alpha| \to \infty$, and
support properties of $\chi_{\alpha},$ we see that $f$ is rapidly decaying along any ray
$r v$ where $v \neq v_\alpha,$ and $f\equiv 0$  within a conic neighborhood of
$v_{\alpha},$ hence  $f\in \cS(\bbR^d).$

\section{Radiation Conditions for Open Wave-Guide Networks}\label{sec.rad_cond}
As noted earlier, the classical Sommerfeld radiation conditions are local on
$\pa\overline{\bbR^d}$: namely, in the absence of channels, $u$ is outgoing in a
neighborhood of $\omega_0 \in S^{d-1}$ if for any sufficiently small
$\delta,\eta>0,$
\begin{equation}
|(\pa_r-ik)u(r\omega)|\leq \frac{C}{r^{\frac{d-1}{2} + \eta}}\text{ uniformly  for }|\omega-\omega_0|\leq \delta.
\end{equation}
In fact, the same condition can be used in the presence of channels so long as
$\omega_0 \not\in \cC$.  In other words, new radiation conditions are only required
near the ends of the channels. As we shall see, these extended radiation
conditions are microlocal in the sense that they require localization in both spatial and frequency variables.

The first successful proposal for radiation conditions for $N$-body
Schr\"odinger operators was given by Isozaki~\cite{Isozaki94}.  His approach
captures the essential ingredients and allows one to prove the key uniqueness
theorem, but is somewhat unwieldy and can be difficult to employ. It also does
not elucidate the fundamental mechanisms that make this the correct
generalization of the Sommerfeld conditions.  A more geometric and flexible
approach was later set forth by Vasy~\cite{VasyAsterisque}, building on and
generalizing the geometric approach pioneered by Melrose~\cite{Melrose94}.  The
assumptions on the potentials used by Isozaki exclude the problem considered in
the present paper, but those of Vasy do cover the case we are considering. In
particular~\cite[condition (11.11)]{VasyAsterisque} allows for potentials of the
form ~\eqref{formq}. In light of this fact, we do not actually use any technical
results from Isozaki's papers, but rather adapt the form of the radiation
condition given in his paper. Its justification relies on Vasy's results.

Isozaki's approach uses more elementary tools from microlocal analysis and is fairly simple to state, but 
only applies to operators on Euclidean space and is overall less flexible. Vasy's approach is based on
Melrose's scattering calculus of pseudodifferential operators and the key idea of the scattering wave-front set, 
which he extends to a more sophisticated and intricate $3$-body scattering calculus with its attendant $3$-body 
scattering wave-front set. This can all be carried out for asymptotically Euclidean (or conical) spaces, 
though we only consider the Euclidean case here. 

We first describe Isozaki's radiation condition, and then introduce the scattering and $3$-body scattering calculi, 
which allow us to then present the Melrose/Vasy formulation.

\subsection{Isozaki's Radiation Conditions}
   In the wave-guide network case the channels are 1-dimensional and meet
   $\pa\overline{\bbR^d}$ at a finite set of points. Because of this we are able
   to first localize to arbitrarily small conic neighborhoods of points on
   $\pa\overline{\bbR^d}$ when checking to see if a formal solution is
   outgoing. This considerably simplifies the needed computations.  In Isozaki's
   work the channels (which are called collision planes) are at least
   2-dimensional, and so meet $\pa\overline{\bbR^d}$ in a positive dimensional
   set; conic localization is therefore not possible near the channel ends.

   We now introduce the class of spatial and frequency localizations that play a
   role in Isozaki's framework.
\begin{definition}
A \emph{conic neighborhood} of $v\in\pa\overline{\bbR^d}$ is a subset of $\overline{\bbR^d}$ of the form
\begin{equation}
\left\{x: \left|\frac{x}{|x|}-v\right|<\delta, \langle x,v\rangle >R\right\},
\end{equation}
with $0<\delta<1,$ and $R > 0.$
\end{definition}
Choose functions $\varphi_{\delta}\in\cC^{\infty}_c(\bbR)$ and  $\psi_{R}\in\cC^{\infty}(\bbR)$ with $\psi_R' \geq 0$ such that
\[
\varphi_{\delta}(t)= \begin{cases}
 &1\text{ for }|t|<\delta,\\
 &0\text{ for }|t|>2\delta,
\end{cases}
\qquad \psi_{R}(t)=
\begin{cases}
&0\text{ for }t<R,\\
&1\text{ for }t>R+1,
\end{cases}
\]
and define the conic cutoff 
\begin{equation}\label{eqn30.203}
\Psi_{v,\delta,R}(x)=\varphi_{\delta}\left(\left|\frac{x}{|x|}-v\right|\right)\psi_{R}(\langle x,v\rangle).
\end{equation}
This has support in a conic neighborhood of $v$. 

\begin{definition}
Fixing $v \in \pa \overline{\bbR^d}$, $\delta, R > 0$ and $u\in\cS'(\bbR^d)$, we say that the tempered distribution
$\cF[\Psi_{v,\delta,R} u]$ is a conic Fourier transform of $u$.
\end{definition}

The distribution $\Psi_{v,\delta,R} u$ is supported in the half space
$\{x_{\alpha}>R\},$ and of polynomial growth, hence for any $s<0,$ the
integral
\begin{equation}\label{eqn31.220}
U_{v_\alpha,\delta,R,s}(\xi)=\int_{\bbR^d}\Psi_{v_\alpha,\delta,R}(x) u(x)e^{-i\langle x,\xi\rangle+s x_{\alpha}}\, dx
\end{equation}
is  absolutely convergent.  Its limit 
\[
U_{v_\alpha,\delta,R}=\lim_{s\to 0^-}
U_{v_\alpha,\delta,R,s}
\]
is well defined as a tempered distribution and equals
$\cF\left[\Psi_{v_\alpha,\delta,R}u\right].$ If $u$ is a formal solution, then is not difficult to see
that $U_{v_\alpha,\delta, R}(\xi)$ is rapidly decreasing as $|\xi|\to\infty,$ though it may not be smooth, or even locally 
represented by a function. For an example see~\eqref{eqn107.201}--\eqref{eqn121.207}.

Isozaki's condition requires one further choice, namely a monotone decreasing
smooth function $\chi_{-,\epsilon}(t)$  that vanishes for $t > 2\epsilon>0,$ and equals $1$ for $t < \epsilon.$ 
  We then say that a solution $u$ is outgoing in a
neighborhood of $v_{\alpha}$ if there exist $\eta, \epsilon, \delta,\delta',R,
R' > 0$ such that, with $\omega=x/|x|,$
\begin{equation}\label{eqn106.211}
\Psi_{v_\alpha,\delta',R'}(x)\cF^{-1}\left[\chi_{-,\epsilon}(\omega\cdot\xi)U_{v_\alpha,\delta,R}(\xi)\right](x)\in  r^{\frac 12-\eta}L^2(\bbR^d).
\end{equation}

This somewhat formidable formula masks a fairly simple idea: namely, we compute
the conic Fourier transform of $u$ around $v_\alpha$, and smoothly cut off the resulting
distribution to a half-space where $\omega\cdot\xi<2\epsilon,$ and then compute the
Fourier inverse.  The solution is outgoing if this new function lies in
$r^{\frac 12-\eta}L^2$ for $\omega$ in a conic neighborhood of
$v_{\alpha}.$ This $L^2$ threshold is based on the fact that $r^{(1-n)/2}$ just
fails to lie in $r^{\frac12} L^2$, so a function which decays at any rate faster
than $r^{(1-n)/2}$ lies not only in $r^{\frac12} L^2$, but in $r^{\frac12 -
  \eta} L^2$ for some sufficiently small $\eta > 0$.

Observe that this condition can also be applied at any point
$v\in\pa\overline{\bbR^d}\setminus\{v_{\alpha}\}.$ At such a point $u$ has an
asymptotic expansion as in Definition~\ref{def1}, hence there exist $\delta, R
>0$ so that $\Psi_{v,\delta,R}(x)u(x)$ has such an expansion. Elementary
calculations then show that $U_{v,\delta,R}(\xi)$ is singular only in a small
neighborhood of $kv\in\bbR^d$; we carry this calculation out in
Appendix~\ref{App1}. Choosing $\epsilon>0,$ small enough, and $\omega$
sufficiently close to $v,$ we can arrange that
$\chi_{-,\epsilon}(\omega\cdot\xi)U_{v,\delta,R}(\xi)\in\cS(\bbR^d).$ From these
observations, \eqref{eqn106.211} follows directly.

We can state Isozaki's outgoing radiation condition as follows:
\begin{definition}[{\bf Isozaki's radiation condition}]\label{def4}
Let $u\in\cS'(\bbR^d)$ be a formal solution to $(H+k^2)u\in\cS(\bbR^d).$ It is defined to be outgoing if, for every $v\in\pa\overline{\bbR^d}$, 
there are positive constants $\epsilon, \delta,\delta',R,R',$ so that 
\begin{equation}\label{eqn108.213}
r^{l} \Psi_{v,\delta',R'}(x)\cF^{-1}\left[\chi_{-,\epsilon}(\omega\cdot\xi)U_{v,\delta,R}(\xi)\right](x)\in L^2(\bbR^d)
\end{equation}
holds for an $l>-\frac 12.$
\end{definition}
\begin{remark}
A solution is \emph{incoming} if it satisfies these conditions, with $\chi_{-,\epsilon}$ replaced by $\chi_{+,\epsilon}(t) :=\chi_{-,\epsilon}(-t).$ 
\end{remark}
\begin{remark}
The condition in~\eqref{eqn108.213} is open in $v,$ hence as
$\pa\overline{\bbR^d}$ is compact, it only needs to be checked at a finite
number of points $\{v_j\}$ for which the interiors of the conic supports of
$\{\Psi_{v_j,\delta'_j,R'}\}$ cover the boundary.
\end{remark}

As we see below, this is essentially a special case of Vasy's condition and
therefore applies to open wave-guide networks.  As in the classical case, the
radiation conditions have many consequences. Vasy shows in~\cite{VasyJFA97} that
a formal solution which is either incoming or outgoing automatically has an
asymptotic expansion
\begin{equation}
u_{\pm}(r\omega)\sim\frac{e^{\pm ikr}}{r^{\frac{d-1}{2}}}\sum_{j=0}^{\infty}\frac{a^{\pm}_j(\omega)}{r^j},
\end{equation}
with smooth coefficients
$\{a^{\pm}_j\}\subset\cC^{\infty}(S^{d-1}\setminus\{v_{\alpha}\}).$ Building on
earlier work of Froese and Herbst, Isozaki, Melrose, et al.\ he proves an
analogue of the Rellich uniqueness theorem, see~\cite{Rellich1943}.
\begin{theorem}[\cite{VasyAsterisque}, Proposition 17.8]\label{thm1}
  If $u\in\cS'(\bbR^d)$ is a solution to $(H+k^2)u=0$ that is outgoing (incoming), then $u\equiv 0.$
\end{theorem}
In Theorem 18.3 of~\cite{VasyAsterisque} Vasy also proves that limits
$$R(k^2\pm i0)=\lim_{s\to 0^+}(H+k^2\pm is)^{-1},$$ 
exist as bounded operators
$$R(k^2\pm i0): H^{0,\delta+\frac{1}{2}}(\bbR^d)\to
H^{0,-\delta-\frac{1}{2}}(\bbR^d),$$ for any $\delta>0.$ We call $R(k^2\pm i0)$
the limiting absorption resolvents. Coupled with the uniqueness theorem above,
we obtain the fundamental existence theorem.
\begin{theorem}[\cite{VasyAsterisque}, Theorem 18.3]\label{thm2}
If $w\in\cS(\bbR^d),$ then the limiting absorption solution $u_+=R(k^2+i0)w$ (resp. $u_-=R(k^2-i0)w$) is the unique outgoing
(resp. incoming) solution to 
\begin{equation}
(H+k^2)u_{\pm}=w.
\end{equation}
\end{theorem}

 The material in Section~\ref{sec4.2}  is rather technical, and
assumes some familiarity with the concepts of microlocal analysis on $\bbR^d.$
A reader can skip to Section~\ref{sec_2d} without much loss of continuity.

\subsection{Vasy's Radiation Conditions}\label{sec4.2}
Isozaki's radiation conditions are specified by classes of
operators\footnote{In~\cite{Isozaki94} a slightly different family of operators
is used, moreover he does not use the conic localization that we employ. See Remark~\ref{rmk7.205}}
$\cR^{\ka}_{\pm}(\epsilon)$ (see~\eqref{eqn108.213}), depending on a set of
parameters, so that a solution to $(H+k^2)u\in\cS$ is outgoing (incoming) if,
for every point $v\in\pa\overline{\bbR^d},$ there is an operator
$P_v\in\cR^{\ka}_{-}(\epsilon),$ (resp.\ $P_v\in\cR^{\ka}_{+}(\epsilon)$ ), `supported at
$v,$' such that $P_vu$ decays faster than $r^{(1-d)/2}$ near $v.$ What is
notable, and somewhat mysterious, about these classes of operators is that they
do not depend explicitly on the free-space wave number $k^2.$ The work of
Melrose and Vasy explains why the conditions introduced by Isozaki provide an
adequate definition for outgoing and incoming solutions. See Remark~\ref{IsozakiExpl}.

Vasy considerably enlarges the class of allowable operators, and this broader
class leads to more refined results.  To state Vasy's radiation conditions we
must describe the 3-body scattering calculus, as defined
in~\cite{VasyAsterisque}, which is a refinement of Melrose's~\cite{Melrose94}
scattering calculus.  We do not describe these calculi in full detail, but give
enough specifics so that the radiation conditions can be accurately stated.
    
\subsubsection{Scattering Calculus} 
The first step is to define the class of scattering differential operators.
These are operators that are naturally defined on any compact
manifold with boundary. We are primarily concerned with the case of operators on
the radial compactification, $\overline{\bbR^d}$, of $\bbR^d,$ which we
identify with the closed unit ball.  This class of differential operators
contains all constant coefficient operators on $\bbR^d,$ among them the
Laplacian and Helmholtz operators.  As well as the Hamiltonians $H =
\Delta_{\bbR^d} + q(x)$  for $q$ a smooth potential, which
decays sufficiently rapidly in all directions.

Consider the radial compactification of Euclidean space. Near the
boundary, which is the `sphere at infinity' $\pa \overline{\bbR^d}$, one may use
`inverted' polar coordinates $\rho=1/|x|,$ and $\omega=x/|x|,$ which we refer
to, in the sequel, as \emph{polar coordinates}.  Near any
$\omega_0\in\pa\overline{\bbR^d},$ there is a subset
$\cJ_{\omega_0}\subset\{1,\dots, d\},$ with $d-1$ elements so that
\begin{equation}
\pa_{x_k}=a_{0k}(\rho,\omega)\rho^2\pa_{\rho}+\sum_{j\in\cJ_{\omega_0}}a_{jk}(\rho,\omega)\rho\pa_{\omega_j}\text{ for }k\in\{1,\dots,d\},
\end{equation}
where the coefficients $\{a_{jk}\}$ are smooth in a neighborhood,
$U_{\omega_0},$ of $\omega_0,$ including at $\rho=0$.
The vector field $\rho^2 \pa_\rho$ has unit length with respect to the Euclidean
metric on $\bbR^d,$ and in this coordinate neighborhood, the vector fields
$\{\rho\pa_{\omega_j}:\:j\in\cJ_{\omega_0}\}$ have norms bounded above and below
by positive constants.  A vector field $X$ on $\overline{\bbR^d}$ belongs to
$\cV_{\Sc}(\overline{\bbR^d}),$ the space of scattering vector fields, if every point
$\omega_0\in\pa\overline{\bbR^d}$ has a neighborhood
$U_{\omega_0}\subset\overline{\bbR^d}$ so that 
\begin{equation}\label{eqn113.217}
X\restrictedto_{U_{\omega_0}}=-\alpha_0(\rho,\omega)\rho^2\pa_{\rho}+\sum_{j\in\cJ_{\omega_0}}\alpha_j(\rho,\omega)\rho\pa_{\omega_j},
\end{equation}
for some smooth functions $\{\alpha_j\}\subset\cC^{\infty}(U_{\omega_0}).$

The duals of these scattering vector fields constitute the space of scattering
$1$-forms; the duals to the specific vector fields above are
$\left\{\frac{d\rho}{\rho^2},\frac{d\omega_j}{\rho}:\:j\in\cJ_{\omega_0}\right\}$.
These are a local spanning set of sections of a new bundle over the compactified
space, $T^*_{\Sc}\overline{\bbR^d},$ called the scattering cotangent bundle.
This bundle agrees with $T^*\bbR^d$ over $\Int \overline{\bbR^d}$; with
\begin{equation}
x=\frac{\omega}{\rho}\Longrightarrow \xi\cdot dx=-\frac{\xi\cdot\omega\, d\rho}{\rho^2}+\frac{\xi\cdot d\omega}{\rho}.
\end{equation}
Carrying this a bit further, note that $|\omega| = 1$ implies $\omega\cdot d\omega=0,$ hence we can use the projections
\begin{equation}
\tau=\xi\cdot\omega,\,\mu=\xi-(\xi\cdot\omega)\omega,
\end{equation}
to define fiber coordinates on $T^*_{\Sc,\omega}\overline{\bbR^d}$ for any $\omega\in U_{\omega_0}.$ A smooth section, $\vartheta$ of
$T^*_{\Sc}\overline{\bbR^d}$ over $U_{\omega_0}$ can be written uniquely as
\begin{equation}
\vartheta=\tau(\rho,\omega)\frac{d\rho}{\rho^2}+\sum_{j\in\cJ_{\omega_0}}\mu_j(\rho,\theta)\frac{d\omega_j}{\rho},
\end{equation}
where $\{\tau,\mu_j\}\subset\cC^{\infty}(U_{\omega_0}).$

We associate to the scattering vector field in~\eqref{eqn113.217} its classical symbol
\begin{equation}
\frac{1}{i} \sigma_1(X)(\rho,\omega;\tau,\mu)= \alpha_0(\rho,\omega)\tau+\sum_{j\in\cJ_{\omega_0}}\alpha_j(\rho,\omega)\mu_j,
\end{equation}
which is a homogeneous polynomial of degree 1 on the fibers of
$T^*_{\Sc}\overline{\bbR^d}.$ This appears to depend on the coordinate choices,
but can be shown to be a coordinate-invariant concept. Over the boundary this
vector field has a \emph{normal symbol}
\begin{equation}
\frac{1}{i}\hN_0(X)(\omega;\tau,\mu)= \alpha_0(0,\omega)\tau+\sum_{j\in\cJ_{\omega_0}}\alpha_j(0,\omega)\mu_j.
\end{equation}
Observe that the radial vector field $\pa_r=-\rho^2\pa_{\rho}$ has normal symbol equal to $i\tau.$

Now define $\Diff^m_{\Sc}(\overline{\bbR^d}),$ the space of scattering
differential operators of order $m,$ as consisting of those operators that can
be expressed as polynomials of order at most $m$ in vector fields $V \in
\cV_{\Sc}(\overline{\bbR^d})$ with coefficients in
$\cC^{\infty}(\overline{\bbR^d}),$ with the convention that
$\Diff^0_{\Sc}(\overline{\bbR^d}) = \cC^{\infty}(\overline{\bbR^d}).$ Just as we
have defined two symbol maps for any $V \in \cV_{\Sc}(\overline{\bbR^d})$, the
same is true for any $P\in\Diff^m_{\Sc}(\overline{\bbR^d})$.  The classical
symbol is the usual homogeneous polynomial of degree $m$ on the fibers of
$T^*\bbR^d$ determined solely by the top order, degree $m$, part of $P$. As
shown above, this extends smoothly to a function on $T^*_{\Sc}\overline{\bbR^d}$
which remains a homogeneous polynomial on the fibers of the scattering cotangent
bundle.

The additional normal symbol, $\hN_0(P),$ is a smooth function on
$T^*_{\Sc}\overline{\bbR^d}\restrictedto_{\pa\overline{\bbR^d}}$.  It is
obtained by replacing $\rho^2\pa_{\rho}$ with $-i\tau$ and $\rho\pa_{\omega_j}$
with $i\mu_j$ and evaluating the coefficients at $\rho=0.$
This is a polynomial of degree $m,$ not generally  homogeneous, on the fibers
of $T^*_{\Sc}\overline{\bbR^d}\restrictedto_{\pa\overline{\bbR^d}}.$ To define
this in a more invariant fashion, choose local coordinates $(\rho,y)$ near any
$p\in\pa\overline{\bbR^d},$ so that the boundary is $\{\rho=0\},$ and $p=(0,0).$
Now choose any smooth function $f$ defined near $p$ and a cutoff
$\psi\in\cC^{\infty}(\overline{\bbR^d}),$ also supported near $p$ with
$\psi(p)=1.$ The limit
\begin{equation}\label{eqn1116.217}
\hN_0(P)(p;\tau,\mu)= \lim_{\rho\to  0^+}e^{-\frac{if(\rho,0)}{\rho}}P[e^{\frac{if(\rho',y')}{\rho'}}\psi](\rho,0),
\end{equation}
is well defined and determined entirely by $\tau=f(0,0)$ and $\mu=d_yf(0,0).$
The two symbols of $P$ are consistent in the sense that the degree
$m$ part of $\hN_0(P)$ equals $\sigma_m(P)\restrictedto_{\pa\overline{\bbR^d}}.$

As a specific example, $P=\Delta+k^2$ has
\begin{equation}
\sigma_2(P)=-(\tau^2+|\mu|^2)\text{ and }\hN_0(P)=k^2-(\tau^2+|\mu|^2).
\end{equation}
Note that $\sigma_2(P)$ is nonvanishing, i.e., elliptic, but the normal symbol vanishes on the set
$$
\Sigma_{k^2}=\{(0,\omega;\tau,\mu):\:k^2=(\tau^2+|\mu|^2)\}\subset
T^*_{\Sc}\overline{\bbR^d}\restrictedto_{\pa\overline{\bbR^d}}.
$$
This is called the \emph{scattering
characteristic variety} of $\Delta+k^2.$ The phenomena of scattering theory
largely take place on this set.

Melrose's next step is to quantize $\Diff^*_{\Sc}(\overline{\bbR^d})$ to a `bi-filtered' algebra of
scattering pseudodifferential operators, $\Psi^{m,l}_{\Sc}(\overline{\bbR^d})$.  The symbols,
$a(x,\xi),$ of operators in $\Psi^{0,0}_{\Sc}(\overline{\bbR^d})$  are elements of $\cC^{\infty}(\overline{\bbR^d}\times
\overline{\bbR^d}),$ i.e., smooth functions on $\bbR^d \times \bbR^d$ which
extend smoothly (and hence are bounded, along with all derivatives) on this
compactification.  This implies the symbolic estimates
\begin{equation}\label{eqnSymbEst}
  |\pa_{x}^{\alpha}\pa_{\xi}^{\beta}a(x,\xi)|\leq \frac{C_{\alpha\beta}}{(1+|x|)^{|\alpha|}(1+|\xi|)^{|\beta|}},
  \text{ for all }\alpha,\beta.
\end{equation}
We identify $\overline{\bbR^d}\times \overline{\bbR^d}$ with
the fiberwise compactification of $T^*_{\Sc}\overline{\bbR^d}.$ An operator
$A\in\Psi^{0,0}_{\Sc}(\overline{\bbR^d})$ acts on $f\in\cS(\bbR^d)$ by the usual
(left quantization) formula
\begin{equation}\label{eqn118.215}
Af(x)=\frac{1}{(2\pi)^d}\int_{\bbR^d}a(x,\xi)e^{ix\cdot\xi}\hf(\xi)d\xi.
\end{equation}

If $r_{\xi}=\sqrt{1+|\xi|^2}$ and $r_x=\sqrt{1+|x|^2},$ then an operator $A\in
\Psi^{m,l}_{\Sc}(\overline{\bbR^d})$ is one that can be represented as
in~\eqref{eqn118.215} with $a\in
r_{\xi}^mr_x^{-l}\cC^{\infty}(\overline{\bbR^d}\times \overline{\bbR^d}).$
Clearly $\Diff^m_{\Sc}(\overline{\bbR^d})\subset
\Psi^{m,0}_{\Sc}(\overline{\bbR^d}).$ We define the symbol $\sigma_m(A)$ for any
such operator $A,$ by extending from $\bbR^d,$ as a homogeneous function of
degree $m$ on the fibers of $T^*_{\Sc}\overline{\bbR^d}$.  Similarly, $\hN_l(A)$
is defined by restricting $r_x^la(x,\xi)$ to $\pa\overline{\bbR^d}.$ It is also
possible to define this normal symbol by oscillatory testing: with $p, (\rho,y), f$
and $\psi$ as in~\eqref{eqn1116.217}, one has
\begin{equation}
\hN_l(A)(p;\tau,\mu)=\lim_{\rho\to  0^+}\rho^{-l}e^{-\frac{if(\rho,0)}{\rho}} A[e^{\frac{if(\rho',y')}{\rho'}}\psi](\rho,0).
\end{equation}
A simple but important observation is that if both $\sigma_m(A)=0$ and
$\hN_l(A)=0,$ then $A\in\Psi_{\Sc}^{m-1,l+1}(\overline{\bbR^d}).$

Define the Sobolev spaces for integral $p\geq 0,$ and all $q\in\bbR$ by
\begin{equation}
H^{p,q}(\bbR^d)=\{u\in\cS'(\bbR^d):\:\pa_x^{\alpha}u\in r_x^{-q}L^2(\bbR^d)\text{ for all }\alpha\text{ with }|\alpha|\leq p\}.
\end{equation}
We can extend the definition to all $p>0$ by interpolation, and to $p<0,$ by
duality: $[H^{p,q}(\bbR^d)]'\simeq H^{-p,-q}(\bbR^d).$ Almost by construction,
any $A\in \Psi^{m,l}_{\Sc}(\overline{\bbR^d})$ defines a bounded map
\begin{equation}
A:H^{p,q}(\bbR^d)\to H^{p-m,q+l}(\bbR^d),\text{ for any }p,q\in\bbR.
\end{equation}
Furthermore, the Fourier transform on tempered distributions restricts to an isomorphism
\begin{equation}
\cF: H^{p,q}(\bbR^d) \longrightarrow H^{q,p}(\bbR^d).
\end{equation}

There is a complete analogue of the familiar symbol calculus in this space of
pseudodifferential operators. An operator $A\in
\Psi^{m,l}_{\Sc}(\overline{\bbR^d})$ is elliptic if both symbols $\sigma_m(A)$
and $\hN_l(A)$ are non-vanishing.  In this case one can follow the usual steps
of the elliptic parametrix construction to obtain an operator
$B\in\Psi^{-m,-l}_{\Sc}(\overline{\bbR^d})$ so that
\begin{equation}
AB-\Id,BA-\Id\in \Psi^{-\infty,\infty}_{\Sc}(\overline{\bbR^d}).
\end{equation}
Note that if $ Q \in \Psi^{-\infty,\infty}_{\Sc}(\overline{\bbR^d})$, then its
Schwartz kernel lies in $\cS(\bbR^d\times\bbR^d).$ More generally, an
operator $A\in \Psi^{m,l}_{\Sc}(\overline{\bbR^d})$ is microlocally elliptic at
$(\tau,\mu)\in T^*_{\Sc\, p}\overline{\bbR^d},$ for $p\in\pa\overline{\bbR^d}$ if
$$
\hN_l(A)(p;\tau,\mu)\neq 0,
$$ and for such an operator one can construct a microlocal parametrix in a
conic neighborhood around $(p; \tau, \mu)$.

Using this calculus of operators we can extend the notion of \emph{wave-front
set}. Traditionally, the wave-front set measures the microlocal {\bf regularity} of a
distribution.  By contrast, the \emph{scattering wave-front set} measures the
\emph{space and frequency localized} {\bf decay} along $\pa\overline{\bbR^d}.$ 

\subsubsection{Classical Wave-front Set}
Before defining the scattering wave-front set, we quickly review the standard
definition of wave-front sets, starting with the even simpler notion of
\emph{singular support}. For any tempered distribution $u\in\cS'(\bbR^d),$ a
point $x_0\notin\singsupp(u)$ if there is a cutoff function
$\psi\in\cC^{\infty}_c(\bbR^d)$ with $\psi \equiv 1$ near $x_0,$ such that $\psi
u\in\cC^{\infty}_{c}(\bbR^d).$ 

If $x_0\in\singsupp(u),$ then $\cF(\psi u)$ is not rapidly decreasing. We can
then assess its behavior in different co-directions.  The set of co-directions
based at a point $x_0$ is identified with the fiber of the co-sphere bundle, $S^*_{x_0}
\bbR^d$, and so the wave-front set, $\WF(u)$, is a subset of the co-sphere
bundle.  Note that $\cF(\psi u)\in\cC^{\infty}(\bbR^d).$ We say that a
co-direction $(x_0,\xi_0)\in S^*_{x_0}\bbR^d$ is \emph{not} in $\WF(u)$ if there
exists a function $\varphi\in\cC^{\infty}(\bbR^d),$ which is homogeneous of
degree 0 outside $B_1(0),$ and equals 1 in a conic neighborhood of the ray
$\{r\xi_0:r\gg 0\},$ for which $\varphi(\xi)\cF(\psi u)(\xi)\in\cS(\bbR^d).$
Thus $\WF(u)$ consists of co-directions $(x_0,\xi_0)$ along which the localized
Fourier transform of $u$ is \emph{not} rapidly decreasing.  This is evidently a
closed set, as its complement is open.

Observe that the composition
\begin{equation}
A_{x_0,\xi_0}: u\mapsto \cF^{-1}[\varphi(\xi)\cF(\psi u)(\xi)]
\end{equation}
is a (classical) pseudodifferential operator of order $0$, with principal symbol
$$
\sigma(A_{x_0,\xi_0})=\varphi(\xi)\psi(x').
$$
This symbol is elliptic in the co-direction $(x_0,\xi_0)$. We thus obtain an alternate definition for the wave-front set:
\begin{definition} A point $(x_0,\xi_0)\notin\WF(u)$ if there exists a (classical) pseudodifferential operator $A\in\Psi^0(\bbR^d),$ with 
$\sigma_0(A)(x_0,\xi_0)\neq 0,$ and  $\psi\in\cC^{\infty}_c(\bbR^d)$ with
  $\psi(x_0)=1,$ such that $A[\psi u]\in\cS(\bbR^d).$ 
\end{definition}

\begin{remark}
  As it can be incorporated into the pseudodifferential operator, $A,$ it is not
  really necessary to include the cut-off function $\psi$ in this definition,
  but we do this here, and in the sequel, to emphasize the local nature of
  wave-front sets.
\end{remark}

We can refine this further and quantify the rate of growth of the localized
Fourier transform in directions belonging to $\WF(u).$ For any $l\in\bbR$, a
co-direction $(x_0,\xi_0)\notin \WF^l(u)$  if there exists a  pseudodifferential operator $A\in\Psi^0(\bbR^d),$ with 
$\sigma_0(A)(x_0,\xi_0)\neq 0,$ and  $\psi\in\cC^{\infty}_c(\bbR^d)$ with
  $\psi(x_0)=1,$ such that $A[\psi u]\in H^l(\bbR^d).$ 
Hence if $(x_0,\xi_0)\in\WF^l(u),$ then, for any operator $A\in\Psi^0$ with
$\sigma_0(A)(x_0;\xi_0)\neq 0,$ we have $A[\psi u]\notin H^l(\bbR^d).$ Each of
these sets is closed, and $\WF^{l'}(u)\subset \WF^{l}(u)$ for any
$l>l'$. Furthermore,
\begin{equation}
\WF(u)=\overline{\bigcup_{l\in\bbR}\WF^l(u)}.
\end{equation}

\subsubsection{Scattering Wave-front Set}
In our setting, $u$ is a \emph{formal solution} of $(H+k^2)u = 0$, that is
$u\in\cS'(\bbR^d),$ and $(H+k^2)u \in\cS(\bbR^d).$ Classical ellipticity implies
that $u\in\cC^{\infty}(\bbR^d)\cap\cS'(\bbR^d)$ and therefore
$\WF(u)=\emptyset.$ To simplify the exposition, we assume this in the sequel.
The `singularities' of interest are something different, namely we are
interested in the directions in $\bbR^d$ along which $u$ does not decay rapidly.
As a simple example, the exponential function $e^{ikx\cdot\omega_0}$ is a smooth
solution to $(\Delta+k^2)u=0,$ but it does not decay at infinity.  The
scattering wave-front set provides a means to localize and quantify this failure to decay.

Proceeding as before, we first define the simpler notion of scattering singular
support. A point $\omega_0\in \pa\overline{\bbR^d}$ is not in the scattering
singular support of $u,$ $\scsingsupp(u),$ if there is a smooth function
$\varphi,$ as above, supported in a conic neighborhood of $\omega_0$  with
$\varphi(\omega_0)=1,$ such that $\varphi u\in\cS(\bbR^d).$ 
We can use the richer structure of
$\Psi^{0,0}_{\Sc}(\overline{\bbR^d})$ to further localize the bad behavior in
frequency as well as spatially along $\pa\overline{\bbR^d}.$
\begin{definition}
A point $(\omega_0;\tau_0,\mu_0)\in T^*_{\Sc}\overline{\bbR^d}\restrictedto_{\pa\overline{\bbR^d}}$ is
\emph{not} in $\WF_{\Sc}(u)$ if there is an operator $A\in\Psi^{0,0}_{\Sc}(\overline{\bbR^d})$ such that 
\begin{equation}\label{eqn125.218}
\hN_0(A)(\omega_0;\tau_0,\mu_0)\neq 0,
\end{equation}
and $A[\varphi u]\in\cS(\bbR^d).$ Here $\varphi$ is smooth,
supported near to $\omega_0,$ with $\varphi(\omega_0)=1.$
\end{definition}
\noindent
If $u$ satisfies $(H+k^2)u=f\in\cS,$ then $u\in H^{m,l}(\bbR^d)$ for all $m\in\bbR,$ and therefore $\cF(u)\in H^{l,m}(\bbR^d)$ 
for all $m \in \bbR$ as well; hence $\cF(u)$ is rapidly decreasing. The scattering wave-front set characterizes
the singularities of $\cF(\varphi u)(\xi).$

A simple choice of operator $A\in\Psi^{0,0}_{\Sc}(\overline{\bbR^d})$ that
satisfies~\eqref{eqn125.218} is defined by
\begin{equation}
Au(z)=\int\tvarphi(z)\psi_0(\omega\cdot\xi)\psi_1(\xi-(\xi\cdot\omega)\omega)\widehat{\varphi u}(\xi)e^{i z\cdot\xi}d\xi,
\end{equation}
where $z=r\omega;$ $\supp\tvarphi\subsubset\supp\varphi\subset B_1(0)^c,$ and
$\psi_0\in\cC^{\infty}_{c}(\bbR),\psi_1\in\cC^{\infty}_{c}(\bbR^d),$ with
$\lim_{r\to\infty}\tvarphi(r\omega_0)=1,$ $\psi_0(\tau_0)=1$ and
$\psi_1(\mu_0)=1.$ This operator has normal symbol
\begin{equation}
\hN_0(A)(z;\tau,\mu)=\tvarphi(z)\psi_0(\tau)\psi_1(\mu).
\end{equation}
As before, there is  a  quantitative refinement.
\begin{definition} Given $m,l\in\bbR,$ and $u\in\cC^{\infty}(\bbR^d)\cap\cS'(\bbR^d),$ a point 
$(\omega_0;\tau_0,\mu_0)\in  T^*_{\Sc}\overline{\bbR^d}\restrictedto_{\pa\overline{\bbR^d}}$ is
\emph{not} in $\WF^{m,l}_{\Sc}(u)$ if there is an operator $A\in\Psi^{0,0}_{\Sc}(\overline{\bbR^d})$ such that 
\begin{equation}
\hN_0(A)(\omega_0;\tau_0,\mu_0)\neq 0,
\end{equation}
and $A[\varphi u]\in H^{m,l}(\bbR^d),$ where $\varphi\in\cC^{\infty}(\overline{\bbR^d}),$ with $\varphi(\omega_0)=1.$
\end{definition}

Melrose introduces the scattering wave-front set in~\cite{Melrose94}, and uses
it define radiation conditions and analyze other properties of solutions to
$(\Delta+q_0(x)+k^2)u\in\cS(\bbR^d),$ where $q_0$ vanishes sufficiently rapidly
at $\pa\overline{\bbR^d}.$ Central to this analysis is the observation that
$\hN_0(\Delta+q_0(x)+k^2)=k^2-\tau^2-|\mu|^2_h$ and therefore $\WF_{\Sc}(u)$ lies
in the \emph{scattering characteristic variety}
\begin{equation}
\Sigma_{k^2}=\{(\omega;\tau,\mu):\: \tau^2+|\mu|_h^2=k^2\}\subset T^*_{\Sc}\overline{\bbR^d}\restrictedto_{\pa\overline{\bbR^d}}.
\end{equation}
 The $\WF_{\Sc}(u)$ is also invariant
under the flow defined on
$T^*_{\Sc}\overline{\bbR^d}\restrictedto_{\pa\overline{\bbR^d}}$ by the rescaled
Hamiltonian vector field defined by $\sigma_2(\Delta)$
\begin{equation}
{}^{\Sc}H_g=2\tau \mu\cdot\pa_{\mu}-2|\mu|_h^2\pa_{\tau}+H_h.
\end{equation}
Here $H_h$ is the Hamiltonian vector field of $h(\omega,\mu),$ the induced metric on $T_{\Sc}^*\pa\overline{\bbR^d}.$ This vector 
field is tangent to $\Sigma_{k^2}$ and vanishes on the two distinguished submanifolds 
\begin{equation}
  R^{\pm}_{k^2}=\{(\omega;\tau,\mu):\:\tau=\pm k,\mu=0\} \subset \Sigma_{k^2},
\end{equation}
which are called the {\it radial sets}.
A trajectory of ${}^{\Sc}H_g$ not contained in $R^+_k\cup R^-_k$ is asymptotic
to a point, $(\omega;k,0),$ on $R^+_k$ as $t\to-\infty,$ and to the 
antipodal point, $(-\omega;-k,0),$ on $R^-_k$ as $t\to \infty.$

Melrose and Vasy (in~\cite{Melrose94,VasyJFA97}) establish the following
properties of $\WF_{\Sc}(u);$  these hold for all real $m:$
\begin{enumerate}
\item If $(\Delta+q_0+k^2)u\in\cS(\bbR^d),$ then $\WF_{\Sc}(u)$ is a union of
  complete trajectories of ${}^{\Sc}H_g.$ If for some $l<-\frac 12, m\in\bbR,$
  $\WF^{m,l}_{\Sc}(u)\cap R^{\pm}_{k^2}\neq\emptyset,$ then this intersection
  must consist of limit points of trajectories of ${}^{\Sc}H_g$ contained
  in $$\WF^{m,l}_{\Sc}(u)\cap[\Sigma_{k^2}\setminus R^+_{k^2}\cup R^-_{k^2}].$$
\item It is possible for $\WF^{m,-\frac 12}_{\Sc}(u)$ to lie entirely in $R^{+}_{k^2}\cup R^-_{k^2}.$ 
\item If for some $l>-\frac 12, m\in\bbR$ $\WF_{\Sc}^{m,l}(u)\cap
  R^{\pm}_{k^2}=\emptyset,$ then $\WF_{\Sc}(u)\cap R^{\pm}_{k^2}=\emptyset.$
\item Finally, if $(\Delta+q_0+k^2)u=0$ and $\WF_{\Sc}^{m,-\frac 12}(u)$ is
  contained in one of $R^{+}_{k^2},$ or $R^-_{k^2},$ then $u\equiv 0.$ This is
  an analogue, and  generalization, of Rellich's uniqueness theorem.
\end{enumerate}

\begin{remark}\label{rmk4.202}
  By adapting the proof of Lemma 11.9 from~\cite{VasyAsterisque}, one can show
  that if $(\Delta+q_0+k^2)u\in\cS(\bbR^d),$ and for some $m,l\in\bbR,$ a
  point $(\omega;\tau,\mu)\notin \WF^{m,l}(u),$ then, for any $m'\in\bbR,$ the
  point  $(\omega;\tau,\mu)\notin \WF^{m',l}(u).$
\end{remark}

Observe that the operators $\pa_r\pm ik$ are elements of $\Psi^{1,0}_{\Sc}(\overline{\bbR^d})$ with 
\begin{equation}
\hN(\pa_r\pm ik)=i(\tau\pm k)
\end{equation}
The classical Sommerfeld radiation condition states that a solution is outgoing
if
\begin{equation}\label{eqn62.202}
(\pa_r- ik)u=O\left(\frac{1}{r^{\frac{d-1}{2} + \delta}}\right)
\end{equation}
for some $\delta > 0$. As $\tau-k\neq 0,$ on $R^{-}_{k^2},$  this implies that $\WF_{\Sc}^{m,l}(u)\cap R^{-}_{k^2}=\emptyset$ for
$l>-\frac 12+\eta,$ for any $\eta<\delta.$  Hence, by 3. above, $\WF_{\Sc}(u)\cap R^{-}_{k^2}=\emptyset.$  
Thus we see that a solution which is ``classically'' outgoing satisfies the condition
\begin{equation}\label{eqn137.218} 
\WF_{\Sc}^{m,l}(u)\cap R^{-}_{k^2}=\emptyset \text{ for an }l>-\frac 12,
\end{equation}
which, in turn, implies that
\begin{equation}
  \WF_{\Sc}(u)\subset R^{+}_{k^2}.
\end{equation}
In~\cite{Melrose94}[\S\S11-12], Melrose also proves that if $u$ satisfies $(\Delta+q_0+k^2)u\in\cS(\bbR^d),$ and
\begin{equation}
  \WF_{\Sc}^{m,-\frac 12}(u)\cap R^{-}_{k^2}=\emptyset,\text{ for some }m\in\bbR,
\end{equation}
then $u$ is classically outgoing and satisfies~\eqref{eqn62.202}.
Similar considerations apply to incoming formal solutions, which satisfy
\begin{equation}
(\pa_r+ ik)u=O\left(\frac{1}{r^{\frac{d-1}{2} + \delta}}\right).
\end{equation}
This is equivalent to
\begin{equation}
  \WF_{\Sc}^{m,-\frac 12}(u)\cap R^{+}_{k^2}=\emptyset,\text{ for some }m\in\bbR,
\end{equation}
and therefore $ \WF_{\Sc}(u)\subset R^{-}_{k^2}.$ This describes the precise
relationship between the classical Sommerfeld radiation conditions and the
scattering wave-front set.

We close this discussion with an example.
\begin{example}
Let us compute the scattering wave-front set of the basic solution
$u_{k,\omega_0}(x)=e^{ik\omega_0\cdot x}$ on $\bbR^2$, where
$\omega_0=(\cos\phi,\sin\phi)\in S^{1}.$ As noted above
$\WF_{\Sc}(u_{k,\omega_0})\subset \Sigma_{k^2}.$ We observe that
\begin{equation}
(\omega_0\cdot\nabla-ik) u_{k,\omega_0}=0.
\end{equation}
The operator $\omega_0\cdot\nabla$ lies in $\Psi^{1,0}_{\Sc}(\overline{\bbR^d})$ and has
\begin{equation}
\hN_0(\omega_0\cdot\nabla)(\theta;\tau,\mu)=i(\tau\cos(\theta-\phi)-\mu\sin(\theta-\phi)).
\end{equation}
The solution to $\hN_0(\omega_0\cdot\nabla)-ik=0$ is the set
\begin{equation}
\begin{split}
\Sigma_{\omega_0} & = \{(\theta;k(\cos(\theta-\phi),-\sin(\theta-\phi)) \\ & 
+l(\sin(\theta-\phi),\cos(\theta-\phi)):\: \theta\in [0,2\pi),l\in\bbR\}.
\end{split}
\end{equation}

From microlocal ellipticity it follows that
\begin{equation}
\begin{split}
\WF_{\Sc}(u_{k,\omega_0})\subset & \Sigma_{k^2}\cap \Sigma_{\omega_0} \\ & = \{(\theta;k(\cos(\theta-\phi),-\sin(\theta-\phi)):\:
\theta\in [0,2\pi)\},
\end{split}
\end{equation}
which consists of two antipodal radial points $(\phi;k,0),(\phi+\pi;-k,0),$ and
the pair of trajectories of ${}^{\Sc}H_g$ joining them.  Note that the
singularities of $e^{ik\omega_0\cdot x}$ over $\phi$ and $\phi+\pi$ are too
strong, according to 1. above, to allow $\WF_{\Sc}(u_{k,\omega_0})$ to be
entirely contained in the radial set, which shows that $\WF_{\Sc}(u_{k,\omega_0})=
\Sigma_{k^2}\cap \Sigma_{\omega_0}.$
Note that the support of the asymptotic expansion of $u_{k,\omega_0}$
in~\eqref{eqn14.210} equals
$\{\pm\omega_0\},$ which is a small subset of $\WF_{\Sc}(u_{k,\omega_0}).$
\end{example}

While the Sommerfeld radiation condition is sufficient to
  define outgoing solutions for a compactly supported (or sufficiently rapidly
  decreasing) potential, such a simple condition cannot suffice in the
  wave-guide network case. The existence of wave-guide modes, shows that there
  are `physical' solutions that do not decay at infinity, and whose radial frequencies
  differs from the free space value. The key
  observation is that the distinction between incoming and outgoing has to do
  only with the sign of the radial frequency ($\tau$) and that there is a critical $L^2$-rate
  of decay. In the remainder of this section we explain the technical framework
  needed to make these notions precise.
  
\subsubsection{3-Body Scattering Calculus}
In order to extend the results above to cover the case when the potential $q$
has channels, Vasy observes that the scattering calculus suffices to define
radiation conditions away from the channel ends, but that the analysis near to
$\cC$ requires the more intricate \emph{3-body scattering calculus}.

The first step in the construction of the 3-body scattering calculus is to pass from the radial compactification of $\bbR^d$ to
a larger compactification where the endpoints of the channels
\[
\cC=\{v_{\alpha}:\: \alpha\in\cA\}\subset\pa\overline{\bbR^d}
\]
are blown up. For each $\alpha\in\cA,$ choose linear coordinates
$(x_\alpha,x^\alpha)$ as in~\eqref{eqn18.222}. For simplicity assume that $v_{\alpha}=(1,0)$ so
that $x_\alpha = x_1$ and the remaining directions are denoted $x'$. With
$(\rho,\omega)$ the corresponding polar coordinates, the channel end
$v_{\alpha}$ is $\{(\rho,\omega):\:\rho=0,\omega_1=1,\omega'=0\}.$ We use $\rho$
and $\omega'$ as coordinates near this point.

To blow up this point we use the polar coordinates centered on $v_{\alpha}:$
\begin{equation}
s=\sqrt{\rho^2+|\omega'|^2},\quad(\zeta_1,\zeta')=\frac{(\rho,\omega')}{s},
\end{equation}
in terms of which the new \emph{front face}, $\ff_{\alpha},$ i.e., the new
boundary component obtained by the blowup, is $\{s=0\}$. The lift of 
$\pa\overline{\bbR^d}\setminus \{v_{\alpha}\}$ is $\{\zeta_1=0\}.$ The new front
face is diffeomorphic to a closed $(d-1)$-dimensional ball. Away from the
$\pa\ff_{\alpha}\approx S^{d-2}$ it is often simpler to use
the coordinates $(1/x_1,x')$ in lieu of $(s, \zeta)$.

\begin{figure}[h]
\centering \includegraphics[width= 12cm]{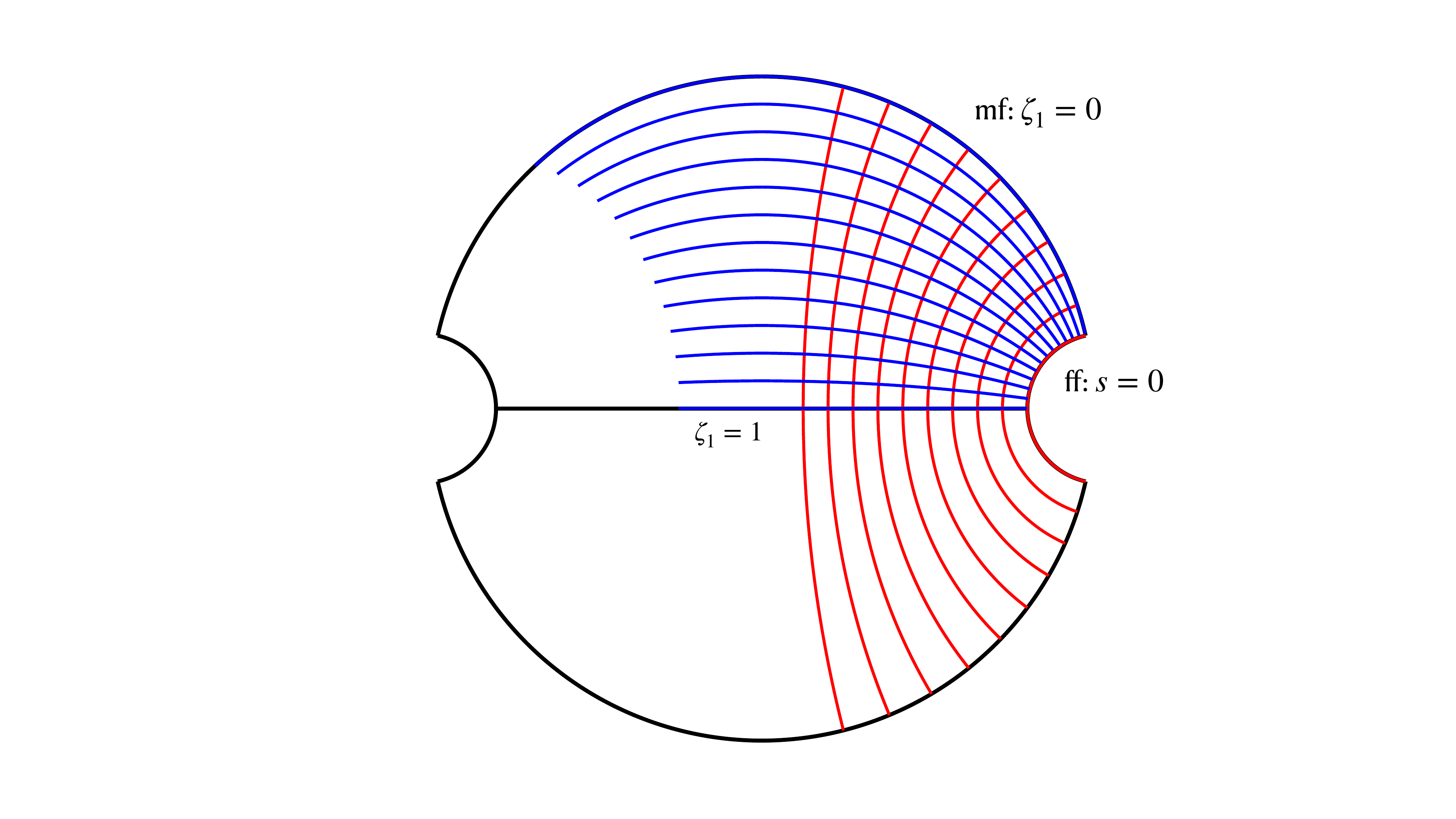}
\caption{Schematic illustration of the wave-guide compactification of
$\bbR^2,$ where $\cC=\{(\pm 1,0)\}.$ Coordinates are shown near the intersection $\mf\cap \ff_{(1,0)}.$ The
red curves are level sets of $s$ and the blue curves are level sets of
$\zeta_1.$ As indicated, in this coordinate neighborhood, $\ff_{(1,0)}=\{s=0\},$ and $\mf=\{\zeta_1=0\}.$ }
\label{figbl_up}
\end{figure}

Following Vasy, we blow up each channel end in this way and denote the new
front face boundary components by $\ff=\cup_{\alpha\in\cA} \ff_{\alpha},$ and
the closure of the lift of $\pa\overline{\bbR^d}\setminus\cC$ by $\mf.$ We call
this new blown-up space the \emph{wave-guide compactification} of $\bbR^d,$ and
denote it by $\overline{\bbR^d}_{\WG}.$ It is equipped with a blow-down map
\begin{equation}
\beta: \overline{\bbR^d}_{\WG}\longrightarrow \overline{\bbR^d}.
\end{equation}
A two-dimensional example appears in Figure~\ref{figbl_up}.

The set of smooth functions on $\overline{\bbR^d}_{\WG}$ is strictly larger than
the space of lifts of elements of $\cC^{\infty}(\overline{\bbR^d}):$ indeed, if
$f$ is smooth on $\overline{\bbR^d},$ then its lift $\beta^*(f)$ is constant on
each $\ff_{\alpha}.$ The potential $q$ in~\eqref{formq} does not extend to be
smooth on $\overline{\bbR^d}$, but it is easy to check that $\beta^*(q)$ does
extend to a smooth function on $\overline{\bbR^d}_{\WG}.$ This is the key
motivation for introducing this blown-up space.  Note also that 
$\beta^*(q_{\alpha})\restrictedto_{\ff_{\alpha}}$ is compactly supported and
$\beta^*(q_{\alpha})\restrictedto_{\ff_{\gamma}}=0$ if $\gamma \neq \alpha$.

Write
$\Diff^0_{3\Sc}(\overline{\bbR^d}_{\WG})=\cC^{\infty}(\overline{\bbR^d}_{\WG}).$
A 3-body scattering vector field is vector field that can be written locally in
the form~\eqref{eqn113.217}, but with coefficients belonging to
$\cC^{\infty}(\overline{\bbR^d}_{\WG}).$ We denote this space of vector fields
by $\cV_{3\Sc}(\overline{\bbR^d}_{\WG})$. An operator $P$ is said to lie in
$\Diff^m_{3\Sc}(\overline{\bbR^d}_{\WG})$ if it can be written locally as a
polynomial of degree $m$ in elements of $\cV_{3\Sc}.$ 
Note that $H=\Delta+q\in\Diff^2_{3\Sc}(\overline{\bbR^d}_{\WG}).$

We now define the symbols for elements of
$\Diff^m_{3\Sc}(\overline{\bbR^d}_{\WG}).$ Both the classical symbol
$\sigma_m(P)$ and, for $p\in\mf,$ the normal symbol $\hN_0(P)(p)$, are defined
as before for $\Diff^m_{\Sc}.$ This part of the normal symbol is a function on
$T^*_{\Sc}\overline{\bbR^d}\restrictedto_{\pa \overline{\bbR^d}\setminus\ff}.$
What remains is to extend the definition of the normal symbol to the front
faces. This is done in Chapter 6 of~\cite{VasyAsterisque}; we summarize
the properties that we need, beginning with the normal symbols of vector
fields. An important property of the normal symbol is that it is multiplicative.

For any $v_{\alpha}\in\cC$, there are coordinates $(\rho,y)$ so that
$\pa\overline{\bbR^d}=\{\rho=0\}$ and $v_{\alpha}$ is the point $(0,0).$ The
variables
\begin{equation}
s=\rho,\quad Y=\frac{y}{\rho},
\end{equation}
are coordinates near $\ff_{\alpha},$ but away from the intersection of  $\ff_{\alpha}$
with $\mf.$ In these blow-up coordinates, the vector fields generating 
$\cV_{\Sc}(\overline{\bbR^d})$ become
\begin{equation}
\rho^2\pa_\rho\rightsquigarrow s^2\pa_s-sY\cdot\pa_Y,\quad \rho\pa_{y_j}\rightsquigarrow \pa_{Y_j}.
\end{equation}
An element $X\in\cV_{3\Sc}(\overline{\bbR^d}_{\WG})$ has the form
\begin{equation}
X= a_0(s,Y)s^2\pa_s+\sum_{j=1}^{d-1}a_j(s,Y)\pa_{Y_j}\text{ with }a_j\in\cC^{\infty}(\overline{\bbR^d}_{\WG}).
\end{equation}
Restricted to $\ff_{\alpha}$, this becomes
$X_{\ff_{\alpha}}=\sum_{j=1}^{d-1}a_j(0,Y)\pa_{Y_j},$ which is a vector field tangent to $\ff_{\alpha}.$ 
This is part of the normal symbol along $\ff_{\alpha}.$  Using the coordinates $(s,\zeta'),$ near $\ff_{\alpha}\cap\mf,$ we see
that it belongs to $\cV_{\Sc}(\ff_{\alpha}).$

If $X,Y\in\cV_{3\Sc}(\overline{\bbR^d}_{\WG}),$ then $[X,Y]\in\cV_{3\Sc}(\overline{\bbR^d}_{\WG})$ and
\begin{equation}\label{eqn77.202}
[X,Y]_{\ff_{\alpha}}=\sum_{j,k} [a_j(0,Y)\pa_{Y_j}b_k(0,Y)-b_j(0,Y)\pa_{Y_j}a_k(0,Y)]\pa_{Y_k}.
\end{equation}
does not vanish on $\ff_{\alpha}.$ We make special note of these important
properties of the normal symbol on a front face: it takes values in a
non-commutative algebra, and the normal symbol of a commutator does not vanish
to one higher order along $\ff_{\alpha}$, as is the case for the normal symbol
along $\mf.$ The normal symbol on $\ff_{\alpha}$ of an element of
$\cV_{3\Sc}(\overline{\bbR^d}_{\WG})$ is a vector field tangent to $\ff_{\alpha},$ which
usually is not constant. In fact, the normal symbol on $\ff_{\alpha}$ should be
understood as a global operator acting on $\cS(\ff_{\alpha}).$

The restriction, $X_{\ff_{\alpha}},$ does not contain any information about the coefficient $a_0(0,Y).$ To recover this information
we use oscillatory testing by considering a function $f\in\cC^{\infty}(\overline{\bbR^d}),$ with $f(0,0)=\tau\in\bbR,\,\pa_y 
f(0,0)=\theta\in\bbR^{d-1}.$ Computing in these coordinates, we see that
\begin{equation}
[e^{-i\frac{f}{s}}Xe^{i\frac{f}{s}}]\psi\restrictedto_{s=0}= e^{-i\theta\cdot Y}\left[X_{\ff_{\alpha}}-i\tau a_0(0,Y)\right]e^{i\theta\cdot Y}\psi(0,Y),
\end{equation}
where $\psi\in\cC^{\infty}(\overline{\bbR^d}_{\WG}).$ There is evidently no loss
of information about the mapping properties of this operator  if we assume
that $\theta=0.$ We define the normal symbol of $X$ on $\ff_{\alpha}$ to be this
family of operators, $\hN_{\ff_{\alpha},0}(X)(\tau)\in
\Psi_{\Sc}^{1,0}(\ff_{\alpha}),$ depending only the $\tau$-variable.  This
variable is essentially a coordinate on the 1-dimensional fiber of the co-normal
bundle to $\ff_{\alpha}.$

This
definition can be extended to $\Diff^m_{3\Sc}(\overline{\bbR^d}_{\WG})$ using
the fact that the symbol is multiplicative,
\[
\hN_{\ff_{\alpha},0}(X\cdot Y)(\tau)=\hN_{\ff_{\alpha},0}(X)(\tau)\cdot\hN_{\ff_{\alpha},0}(Y)(\tau).
\]
For example, in terms of these projective coordinates,
$\hN_0(\pa_{x_1})(\tau)=i\tau,$ and
\begin{equation}
  \hN_{\ff_{\alpha},0}(H+k^2)(\tau)=\Delta_Y+q_{\alpha}(Y)+k^2-\tau^2.
\end{equation}

In~\cite{VasyAsterisque} Vasy introduces a bi-graded algebra of
pseudodifferential operators, $\Psi^{m,l}_{3\Sc}(\overline{\bbR^d}_{\WG}),$
which quantizes $\Diff^*_{3\Sc}(\overline{\bbR^d}_{\WG})$. An operator
$A\in\Psi^{m,l}_{3\Sc}(\overline{\bbR^d}_{\WG})$ has both a classical symbol,
$\sigma_m(A),$ and a normal symbol $\hN_l(A).$ The normal symbol agrees with the
previous definition over $\Int\mf$ and is a function on the fibers of
$T^*_{\Sc}\overline{\bbR^d}\restrictedto_{\pa\overline{\bbR^d}\setminus\ff}.$

Over a front face $\ff_{\alpha}$, the normal symbol $\hN_l(A)(v_{\alpha};\tau)$ is a
function of $\tau\in\bbR$ that takes values in $\Psi^{m,0}_{\Sc}(\ff_{\alpha}).$
`Nonvanishing' of the symbol at $\tau_0$ is  interpreted to mean that
the operator $\hN_l(A)(v_{\alpha};\tau_0),$ acting on $L^2(\ff_{\alpha}),$ is invertible. In
this case the inverse belongs to $\Psi^{-m,0}_{\Sc}(\ff_{\alpha}).$  The
fact that we need to define these symbols as operators acting on
$L^2(\ff_{\alpha})$ is a reflection of the
fact, already evident from~\eqref{eqn77.202},  that the normal symbol map over $\ff_{\alpha}$ takes
values in a non-commutative algebra.

We now extend the definition of the wave-front set to the front face.
\begin{definition} For $\alpha\in\cA,$
a point $\tau_0\in\bbR$ is not in $\WF^{m,l}_{3\Sc, \ff_{\alpha}}(u),$ if there
is an operator $A\in\Psi_{3\Sc}^{0,0},$ a function
$\varphi\in\cC^{\infty}(\overline{\bbR^d})$ with $\varphi(v_{\alpha})=1,$ for
which $\hN_0(A)(v_{\alpha};\tau_0)$ invertible and
\begin{equation}
A(\varphi u)\in H^{m,l}(\bbR^d).
\end{equation}
\end{definition}

The normal symbol of $H+k^2$ over $\mf$ is $\hN_{\mf,0}(H+k^2)=k^2-(\tau^2+|\mu|^2);$ the characteristic variety over 
$\mf$ is the closure, in $\pa(\overline{\bbR^d}_{\WG}),$ of
\begin{equation}
\Sigma_{k^2,\mf}=\{(\omega;\tau,\mu):\:\omega\notin\beta^{-1}(\cC);\,\tau^2+|\mu|^2=k^2\}.
\end{equation}
The normal symbol of $H+k^2$ over $\ff_{\alpha}$ equals
\begin{equation}
\hN_{0}(H+k^2)(v_{\alpha};\tau)=\Delta_{\ff_{\alpha}}+q_{\alpha}(Y)+k^2-\tau^2,
\end{equation}
where $\Delta_{\ff_{\alpha}}$ is the Laplacian  with respect to the Euclidean
metric defined in $\bbR^{d-1}.$  This allows us to extend the characteristic variety over the front faces. On $\ff_{\alpha}$ it is 
the set of $\tau\in\bbR$ for which $\hN_{0}(H+k^2)(v_{\alpha};\tau)$ is not invertible.

It is well known that $\hN_{0}(H+k^2)(v_{\alpha};\tau)$ is invertible except where $|\tau|\leq k$ or
$\tau\in\left\{\pm\sqrt{k^2+E^2_{\alpha,j}}:\: j=1,\dots, N'_{\alpha}\right\}.$ Proposition 11.2 in~\cite{VasyAsterisque} 
states that the characteristic variety of $H+k^2$ over $\ff_{\alpha}$ is the union of the two sets:
\begin{equation}\label{eqn80.201}
\Sigma_{k^2,\ff_{\alpha}}^{\pm}=\left\{0\leq\pm\tau\leq k;\,\tau=\pm\sqrt{k^2+E^2_{\alpha,j}}:\;j=1,\dots, N'_{\alpha}\right\}.
\end{equation}
As shown in~\cite{VasyAsterisque}, if $(H+k^2)u\in\cS,$ then $\WF_{3\Sc}(u)\restrictedto_{\ff_{\alpha}}$ is contained in 
the union of the these two sets.

If a  formal solution, $u\in\cS'(\bbR^d)$ with $(H+k^2)u\in\cS,$ 
satisfies a standard radiation condition away from the channels ends, then it
follows from  Proposition 14.1 of~\cite{VasyAsterisque} that 
\begin{equation}\label{eqn84.202}
 \WF_{3\Sc}(u)\restrictedto_{\ff_{\alpha}}\subset
   \bigcup_{j=1,\dots,
    N'_{\alpha}}\left\{\pm \tau=\sqrt{k^2+E^2_{\alpha,j}}\right\}\cup  \{\pm \tau = k\}.
\end{equation}
Indeed Vasy refined the propagation results of Melrose to show that if a point $(v_{\alpha};\tau),$
with $|\tau|<k$ belongs to $\WF_{3\Sc}(u)\restrictedto_{\ff_{\alpha}},$ then it must be the limit point of
 a   trajectory of ${}^{\Sc}H_g$ contained in $\WF_{3\Sc,\mf}(u)\setminus R^+_{k^2}\cup R^-_{k^2}\restrictedto_{\mf}.$ On the other
 hand, the existence of such a trajectory is impossible if
 $$\WF_{3\Sc,\mf}(u)\subset
 R^+_{k^2}\cup R^-_{k^2}\restrictedto_{\mf}.$$
 The radial sets over the front faces $\{\ff_{\alpha}:\:\alpha\in\cA\}$ are
 defined to be
\begin{equation}\label{eqn85.203}
R^{\pm}_{k^2,\ff_{\alpha}}= \left\{\pm\tau=k;\, \tau=\pm\sqrt{k^2+E^2_{\alpha,j}}:\;j=1,\dots, N'_{\alpha}\right\}.
\end{equation}

The definition of $\WF_{3\Sc}(u)$ is microlocal. A formal solution to
$(H+k^2)u\in\cS$ is outgoing over $\Int\mf$ if it satisfies the condition
in~\eqref{eqn137.218} restricted to points lying over $\Int\mf.$ For a point
$(\omega_0,\xi_0)\in [R^+_{k^2}]^c,$ with $\omega_0\in\Int\mf,$ this means that
there is an operator $A\in\Psi^{*,0}_{\Sc}(\overline{\bbR^d})$ for which
$\hN_0(A)$ is supported in $\Int\mf,$ and $\hN_0(A)(\omega_0;\xi_0)\neq 0$ so
that, for some $l>-\frac12,$ and $\varphi\in\cC^{\infty}(\overline{\bbR^d})$
with $\varphi(\omega_0)=1,$ we have
\begin{equation}
A(\varphi u)\in H^{*,l}(\bbR^d).
\end{equation}
Thus we only need to give a outgoing radiation condition for points lying above the
channel ends. From~\eqref{eqn84.202} we see that we only need to show that
  the $\WF_{3\Sc,\ff}(u)$ lies in  $\{\tau>0\}.$

\begin{remark}\label{IsozakiExpl}
  These observations contain the explanation for the adequacy of Isozaki's
  conditions: the radial sets lie in disjoint half spaces defined by
  $\{\tau<0\},$ and $\{\tau>0\},$ where $\tau=x\cdot\xi/|x|.$ If there were
  points in $\WF_{3\Sc}(u)$ not in the radial set, then the propagation results
  proved in~\cite{Melrose94,VasyAsterisque} would imply that a non-trivial,
  complete (broken) trajectory of ${}^{\Sc}H_g$ would also have to belong to
  $\WF_{3\Sc}(u).$ But no such trajectory is contained in a half space
  $\{\tau<0\},$ or $\{\tau>0\}.$ Hence the only way that $\WF_{3\Sc}(u)$ can be
  contained in such a half space is if it is contained in a radial set.
    Furthermore, there is a decay threshold that is required in order for these
    wave-fronts to be contained in $R^+_{k^2}\cup R^-_{k^2},$ which explains the
    requirement that $l>-\frac 12$ in~\eqref{eqn108.213}. On the $\mf$ the existence of this threshold is the content of 
    point 1. preceding Remark~\ref{rmk4.202}. 
\end{remark}

\begin{definition}[{\bf Vasy's Radiation Condition}]
A formal solution $u\in\cS'$ of $(H+k^2)u\in\cS$ is outgoing if it is outgoing over $\Int\mf$ and, for every $v_{\alpha}\in\cC$ and 
$\tau_0\in\Sigma_{k^2,\ff_{\alpha}}^{-}$, there is an operator
$A\in\Psi^{0,0}_{3\Sc}$ with $\hN_0(A)(v_{\alpha};\tau_0)$ invertible, a function
$\varphi\in\cC^{\infty}(\overline{\bbR^d})$ with $\varphi(v_{\alpha})=1,$  such
that, for some $l>-\frac 12,$ and $m\in\bbR,$
\begin{equation}\label{eqn85.220}
A(\varphi u)\in H^{m,l}(\bbR^d).
\end{equation}
\end{definition}
\begin{remark}
Since formal solutions belong to $H^{m,l},$ for all $m\in\bbR,$ we can use
operators $A\in\Psi^{m,0}_{3\Sc}(\overline{\bbR^d}_{\WG})$ with
$\hN_0(A)(v_{\alpha};\tau_0)$ invertible, for any $m\in\bbR.$ The symbol of the
operator $\hN_0(A)(v_{\alpha};\tau_0)\in\Psi^{0,0}_{\Sc}(\ff_{\alpha})$ extends
continuously to $\pa\ff_{\alpha};$ the invertibility of this normal operator implies
that $\hN_0(A)(\omega;\tau,\mu)$ is non-vanishing for
$\omega\in\pa\ff_{\alpha},$ with $|\tau-\tau_0|,$ and $|\mu|$ sufficiently small.
\end{remark}
\begin{remark}\label{rmk7.205}
  The operators, $A_{{\pm}},$ appearing in Isozaki's radiation conditions have (classical) ``triple-symbols''
\begin{equation}\label{eqn74.202}
a_{{\pm}}(x,\xi,x')=  \tvarphi(x)\chi_{\pm}\left(\frac{x\cdot \xi}{|x|}\right)\varphi(x'),
\end{equation}
where $\chi_\pm(t)$ vanishes for $\mp t > 2\epsilon$ and equals $1$ for $\mp t <
\epsilon$.  The functions $\varphi,\,\tvarphi\in\cC^{\infty}(\overline{\bbR^d})$
are conic cutoffs, with $\tvarphi(x)=\varphi(x)=1$ for $x$ near $v_{\alpha},$
with supports contained in $B_1^c(0),$ see~\eqref{eqn30.203}. In this case
\begin{equation}
\hN_{\ff_{\alpha},0}(A_{\pm})(\tau)=\chi_{\pm}(\tau).
\end{equation}

Observe that these operators do not belong to
$\Psi^{0,0}_{\Sc}(\overline{\bbR^d}),$ (or
$\Psi^{0,0}_{3\Sc}(\overline{\bbR^d}_{WG})$) because $x$-derivatives lead to increased
growth in $\xi$.  
  In fact, in~\cite{Isozaki94} Isozaki defines operator/symbol classes, $\cR^{\ka}_{\pm}(\epsilon),$ as follows: a function
  $p(x,\xi)\in \cC^{\infty}(T^*\bbR^d)$ belongs to $\cR^{\ka}_{\pm}(\epsilon),$  if it satisfies finitely many symbolic estimates
  \begin{equation}
    |\pa_x^{\alpha}\pa_{\xi}^{\beta}p(x,\xi)|\leq\frac{C_{mn}}{(1+|x|)^{|\alpha|}(1+|\xi|)^{\ka}},\,\text{ for }0\leq |\alpha|, |\beta|\leq \ka,
  \end{equation}
  and, on $\supp p,$
    \begin{equation}
      \inf_{x,\xi}\pm \frac{x\cdot \xi}{(1+|x|)}> -\epsilon.
    \end{equation}
    The associated operator is defined as usual by
    \begin{equation}
      Pv(x)=\cF^{-1}\left[ p(x,\xi)\hv(\xi)\right](x).
    \end{equation}

As observed by Vasy, a formal solution $u\in H^{m,l}$ for any $m$, and therefore
$\widehat{\varphi u} \in H^{l,m}.$ Hence, for any $n\in\bbN,$  we can  rewrite
\begin{equation}
A_{\pm}u(x)=\tvarphi(x)\cF^{-1}\left[\frac{\chi_{\pm}\left(\frac{x}{|x|} \cdot\xi \right)}{(1+|\xi|^2)^n}\cdot (1+|\xi|^2)^n
\widehat{\varphi u}(\xi)\right].
\end{equation}
For any $p\geq 0,$ by choosing $n$ large enough, we can arrange to have the symbols
\begin{equation}\label{eqn149.218}
\frac{\chi_{\pm}\left(\frac{x}{|x|} \cdot\xi \right)\varphi(x)}{(1+|\xi|^2)^n},
\end{equation}
satisfy any finite number of symbolic estimates for an operator in
$\Psi^{-p,0}_{\Sc}(\overline{\bbR^d}),$ see~\eqref{eqnSymbEst}. It is in this sense that Isozaki's
operator classes $\cR^{\ka}_{\pm}(\epsilon)$ belong to
$\Psi^{*,0}_{\Sc}(\overline{\bbR^d}),$ and are a subset of the operators used
by Vasy.
\end{remark}
\begin{remark}
Vasy's radiation condition is very much like Isozaki's except that we can test
to see if a formal solution is outgoing with operators
$A\in\Psi^{m,0}_{3\Sc}(\overline{\bbR^d}_{WG}).$ This is a much larger and more
flexible class than the classes $\cR^{\ka}_{\pm}(\epsilon).$ From Vasy's
  formulation we see that to decide if a point $(\omega;\xi)\in
  T^*_{\Sc}\overline{\bbR^d}_{\WG}\restrictedto_{\overline{\bbR^d}_{\WG}}$
  belongs to $\WF_{3\Sc}^{m,l}(u),$ what is needed is an operator $A$ whose
  normal symbol, $\hN_0(A)$ is invertible at this point.   An interesting
example of a differential operator, useful for this purpose, which belongs to
$\Psi^{1,0}_{3\Sc}(\overline{\bbR^d}_{\WG})\setminus
\Psi^{1,0}_{\Sc}(\overline{\bbR^d})$ is given in~\eqref{eqn198.219}.  Once we
know that $\WF_{3\Sc}(u)\restrictedto_{\Int\mf}\subset [R^+_{k^2}\cup
  R^-_{k^2}]_{\Int\mf},$ then, near to $\ff_{\alpha},$ it suffices to use
operators like $A_{-}$ defined by~\eqref{eqn74.202}, where for any $\epsilon>0,$
we take $\chi_-\in\cC_c^{\infty}(\bbR),$ with
\begin{equation}
\chi_-(t)=
\begin{cases}
&0\text{ for }t\notin[-\sqrt{k^2+E^2_{\alpha,N'_{\alpha}}}-2\epsilon,-k+2\epsilon]\\
&1\text{ for }t\in [-\sqrt{k^2+E^2_{\alpha,N'_{\alpha}}}-\epsilon,-k+\epsilon].
\end{cases}
\end{equation}
While these operators still do not belong to
$\Psi^{*,*}_{3\Sc}(\overline{\bbR^d}_{\WG}),$ as Vasy explains in Chapter 13
of~\cite{VasyAsterisque},  they can be treated as if they do. 
\end{remark}

As noted earlier, Vasy shows~\cite[Theorem 18.3]{VasyAsterisque} that the limiting absorption resolvents exist as maps
\[R(k^2\pm i0)=\lim_{\sigma\to 0^+} R(k^2\pm i\sigma): H^{0,\frac 12+\delta} \longrightarrow H^{0,-\frac 12-\delta}
\]
for any $\delta>0.$ Moreover for any $f\in\cS,$ 
\begin{equation}\label{eqn85.201}
  \WF_{3\Sc,\mf}[R(k^2\pm i0)f]\subset R^{\pm}_{k^2}\restrictedto_{\mf}\text{ and }
  \WF_{3\Sc,\ff_{\alpha}}[R(k^2\pm i0)f]\subset R^{\pm}_{k^2,\ff_{\alpha}}.
\end{equation}
\begin{remark}
Note that the signs of $i0$ are switched relative to those used
in~\cite{VasyAsterisque} since we work with $H=\Delta+q$ rather 
than $H=-\Delta+q.$
\end{remark}
 
\section{The Integral Equation Solutions to the $2d$--Model Problem}\label{sec_2d}
In this section we show that the solutions to the 2d-model problems defined
in~\cite{EpWG2023_1} and analyzed in~\cite{EpWG2023_2} satisfy the outgoing
radiation conditions for an open wave-guide network, as described in
Section~\ref{sec.rad_cond}. We make extensive use here of the notation,
definitions, and results from~\cite{EpWG2023_1} and~\cite{EpWG2023_2}, where it
is assumed that there are two channels, both co-linear with the $x_1$ axis. In
these papers, the potential $q(x)$ is assumed to be piecewise constant, equal to
a `left' potential $q_l(x_2)$ when $x_1 < 0,$ and a `right' potential $q_r(x_2)$
when $x_1 > 0$, with a jump discontinuity contained in a finite interval,
$\{0\}\times[-d,d],$ of the $x_2$-axis,  that is
\begin{equation}\label{eqn98.207}
  q(x_1,x_2)=q_l(x_2)\chi_{(-\infty,0]}(x_1)+q_r(x_2)\chi_{[0,\infty)}(x_1).
\end{equation}
See Figure~\ref{fig0}. In this paper we need to assume that $q_l$ and $q_r$ are
both smooth and compactly supported, though $q,$ defined in~\eqref{eqn98.207},
still has a jump discontinuity along a finite interval of the $x_2$-axis.

The arguments in~\cite{EpWG2023_1,EpWG2023_2}, readily extend to handle such `piecewise smooth'
potentials, although the formulas become slightly less explicit, see~\cite{GHRQ}. The results of
Section~\ref{sec.rad_cond} herein are unaffected by the presence of a compactly
supported jump discontinuity in the potential. It is very likely that the radiation
conditions and uniqueness results extend to piecewise constant, or even more
singular potentials.  The microlocal proofs would need to be replaced by
something similar to that used in the classical $N$-body scattering literature,
where singular potentials are routinely considered, see~\cite{PSS1981}.

As before, we  assume that $v\equiv 0$ is the only bounded solution, on $\bbR,$ of either of
the two equations $(\pa_{x_2}^2+q_{l,r}(x_2))v(x_2)=0.$ That is, there is no
threshold at $k=0.$ This technical hypothesis is needed to be able to obtain the
representations for the limiting absorption resolvents $(\Delta+q_{l,r}+k^2+i0)^{-1}$ as
contour integrals, which are used extensively in~\cite{EpWG2023_1}
and~\cite{EpWG2023_2}.  The radiation conditions can be phrased in terms of
conically localized Fourier transforms, as in Isozaki~\cite{Isozaki94}, or using
the characterization of the 3-body scattering wave-front set, as
in~\cite{VasyAsterisque}.  We make use of both formulations in the computations that follow.

\subsection{Scattering as a Transmission Problem}
In~\cite{EpWG2023_1,EpWG2023_2}, the problem of scattering an incoming wave
  off of an open wave-guide network is rephrased as a transmission problem.
  We now briefly outline this approach.  Assume that we are given left and
right incoming fields, $u^{\In}_{l,r}$ that satisfy the equations
$$(\Delta+k^2+q_{l,r}(x_2))u^{\In}_{l,r}=0,\text{ where }\pm x_1\geq 0.$$
The scattered fields are then found as
solutions to the equations
\begin{equation}\label{eqn99.221}
  (\Delta+k^2+q_{l,r}(x_2))u^{l,r}=0,\text{ where }\pm x_1\geq 0,
  \end{equation}
which satisfy the transmission boundary conditions along the $x_2$-axis:
\begin{equation}\label{eqn88.220}
  \begin{split}
   &g(x_2) \overset{d}{=}
    u^r(0^+,x_2)-u^l(0^-,x_2)=u^{\In}_{l}(0^+,x_2)-u^{\In}_r(0^-,x_2)\\
&h(x_2)\overset{d}{=}
    \pa_{x_1}u^r(0^+,x_2)-\pa_{x_1}u^l(0^-,x_2)=\pa_{x_1}u^{\In}_{l}(0^+,x_2)-\pa_{x_1}u^{\In}_r(0^-,x_2),
    \end{split}
\end{equation}
and are ``outgoing'' in some sense. Below we show that the scattered fields we obtain
  are, in fact,  outgoing in the sense defined in Section~\ref{sec.rad_cond}.

This transmission problem is solved by
  first finding the kernels, $\cE^{l,r}(x;y),$ for outgoing fundamental
  solutions, $(\Delta+q_{l,r}+k^2+i0)^{-1},$ for the
  bi-infinite wave-guide operators $(\Delta+q^{l,r}(x_2)+k^2).$  
   These kernels
  are given by  fairly explicit formul{\ae}, as perturbations of the free-space
  fundamental solution
  \begin{equation}
    \cE^{l,r}(x;y)=g_k(x-y)+\cW^{l,r}(x;y),
  \end{equation}
  where
  \begin{equation}
    g_k(x-y)=\frac{i}{4\pi}H^{(1)}_0(k|x-y|),
  \end{equation}
  is the free space kernel. We let $\cS_k$ and $\cD_k$ denote the classical single and
  double layer operators over the $x_2$-axis defined by $g_k.$

  In each half plane, $\pm x_1>0,$ we represent
  \begin{equation}\label{eqn103.700}
    u^{l,r}(x_1,x_2)=\int_{-\infty}^{\infty}\left[\cE^{l,r}(x;0,y_2)\tau(y_2)-\pa_{y_1}\cE^{l,r}(x;0,y_2)\sigma(y_2)\right]dy_2.
  \end{equation}
  Here and in the sequel $\sigma$ and $\tau$ are functions defined on the real
  line with adequate decay for these integrals to be absolutely
  convergent. Letting $x_1\to 0^{\pm}$ in these representations, and using the
  transmission conditions,~\eqref{eqn88.220}, we obtain integral equations for
  $(\sigma,\tau),$ see~\eqref{eqn122.203}. Modulo the uniqueness theorem proved
  below (Theorem~\ref{thm3}), we show in~\cite{EpWG2023_1} that these integral
  equations are solvable for physically interesting incoming data.

  These solutions are represented in the appropriate half planes as sums of
  terms corresponding to the decompositions of $\cE^{l,r}$ described above:
  \begin{equation}\label{eqn89.220}
    u^{l,r}(x_1,x_2)=[u^{l,r}_0-u^{l,r}_1+u^{l,r}_{c0}-u^{l,r}_{c1}+u^{l,r}_g](x_1,x_2)\text{
      where }\pm x_1\geq 0.
  \end{equation}
  Here $u^{l,r}_0-u^{l,r}_1=\cS_k\tau-\cD_k\sigma$ are the contributions of the
  `free space' Green's functions,
  $u^{l,r}_{c0}-u^{l,r}_{c1}=\cW^{l,r}_c\tau-\cW^{l,r\, '}_c\sigma,$ are the
  ``continuous spectral'' parts of the correction terms needed to obtain the kernels for
  $(\Delta+q_{l,r}+k^2+i0)^{-1},$ and $u^{l,r}_g$ are the contributions of the
  wave-guide modes.

  The  function
  \begin{equation}\label{eqn90.220}
    u^{\tot}=
    \begin{cases}
      &u^{l}(x_1,x_2)+u^{\In}_l(x_1,x_2)\text{ where }x_1\leq 0,\\
       &u^{r}(x_1,x_2)+u^{\In}_r(x_1,x_2)\text{ where }x_1\geq 0,
    \end{cases}
  \end{equation}
  is a weak solution to $(\Delta+q(x_1,x_2)+k^2)u^{\tot}=0,$  where
\begin{equation}
q(x_1,x_2)=\chi_{(-\infty,0]}(x_1)q_l(x_2)+\chi_{[0,\infty)}(x_1)q_r(x_2).
\end{equation}
It has, at worst, a jump discontinuity in the $\pa^2_{x_1}$-derivative along a
finite interval of the $x_2$-axis. 

We now show that the various parts of $u^l$ and $u^r$ are outgoing. We begin
with the wave-guide modes, i.e., a solution of the form
$u(x_1,x_2)=e^{i\xi^\circ x_1}v(x_2),$ to either the left or right equations
\begin{equation}
(\Delta+q_{l,r}(x_2)+k^2)u=0.
\end{equation}
Here $\xi^\circ = \pm \sqrt{k^2 + E^2}$, and $(\pa_{x_2}^2 + q_{l,r}(x_2))
v(x_2) = E^2 v(x_2)$.  Below we show that, if $\xi^\circ>0,$ then this solution
is outgoing to the right, while if $\xi^\circ < 0$, then it is outgoing to the
left. Since $v(x_2)$ is smooth and exponentially decaying, a conical Fourier
transform localized to any non--horizontal ray belongs to Schwartz class. Hence
the scattering wave front set of $u$ can lie only over the points $(\pm
1,0)\in\pa\overline{\bbR^2}.$ Denote the front faces of the corresponding
blowups by $-\leftrightarrow \ff_l$ and $+\leftrightarrow \ff_r$.

In Vasy's approach, the normal symbol of the operator
$\pa_{x_1}-i\xi^\circ\in\Psi^{1,0}_{3\Sc}$ on $\ff_r$ is $i(\tau-\xi^\circ)$ and
on $\ff_l$ it is $-i(\tau+\xi^\circ)$. These symbols  vanish on $\ff$ at  the points $\{
(1,0;\xi^\circ), (-1,0;-\xi^\circ)\}$.  Since $(\pa_{x_1}-i\xi^\circ) u=0,$ it
follows from~\eqref{eqn80.201} that
\begin{equation}
\WF_{3\Sc}(u)\subset \{(-1,0;-\xi^\circ), (1,0;\xi^\circ)\}.
\end{equation}
Hence, if $\xi^\circ>0,$ then  $u$ is outgoing as $x_1\to\infty$ and incoming as
$x_1\to-\infty,$ and vice-versa if $\xi^{\circ}<0.$

Applying Isozaki's formulation takes a bit more work. We carry out the analysis
as $x_1\to\infty$. Choose $\psi_+\in\cC^{\infty}(\bbR)$ monotone with $\psi_+(x)
= 0$ for $x \leq -1$ and $\psi_+(x) = 1$ for $x \geq 1$, and assume that
$\pa_{x_1}\psi$ is even.  This is a conic cutoff, as defined in~\eqref{eqn30.203}, where
the angular term is identically 1.  The Fourier transform of
$\psi_+(x_1)u(x_1,x_2)$ is well defined as a distribution, with an analytic
continuation to the lower half plane in the $\xi_1$-variable. At frequency $(\xi_1+is),$ with $s<0,$  the Fourier integral is absolutely
convergent  and equals
\begin{equation}\label{eqn107.201}
\begin{split}
\widehat{\psi_+u}(\xi_1+is,\xi_2) & =\int_{\bbR^2} \psi_+(x_1)e^{-i[\xi_1-\xi^\circ+is]x_1}v(x_2)e^{-i\xi_2 x_2}\, dx_1dx_2 \\
& = \widehat{\psi_+}(\xi_1 - \xi^\circ + is) \hat{v}(\xi_2).
\end{split}
\end{equation}
This in turn converges, distributionally, as $s \to 0^-$ to $\widehat{\psi_+}(\xi_1 - \xi^\circ) \hat{v}(\xi_2)$.

It is useful to have a more explicit description of the singularity of
$\widehat{\psi_+}$.  Returning to its regularized form, with $s < 0$, we
integrate by parts to obtain
\begin{equation}\label{eqn105.204}
\begin{split}
\widehat{\psi_+ }(\xi_1+is) & =\int_{-\infty}^{\infty} \psi_+(x_1)e^{(s-i\xi_1)x_1}\, dx_1 \\ & =-\int_{-\infty}^{\infty} \frac{\psi'_+(x_1)e^{(s-i\xi_1)x_1}}{i(\xi_1-is)}\, dx_1= -\frac{\widehat{\psi'_+}(\xi_1+is)}{i(\xi_1-is)}.
\end{split}
\end{equation}
Since $\psi_+'(x_1)$ is compactly supported and even,
$\widehat{\psi'_+}(\xi)\in\cS(\bbR)$ is real-valued and also even, with
$\widehat{\psi'_+}(0)=1.$ Thus if $\varphi \in\cC^{\infty}_c(-2\delta,2\delta)$
is even, with $\varphi(t)=1$ for $|t|<\delta,$ then, for $s\leq 0,$
\begin{equation}
\widehat{\psi_+}(\xi_1+is) + \frac{\widehat{\psi'_+}(\xi_1+is)\varphi(\xi_1)}{i(\xi_1-is)} \in \cS(\bbR),
\end{equation}
and its rapid decrease is uniform as $s\to 0^-$. As distributions, 
 \begin{equation}\label{eqn121.207}
\lim_{s\to 0^-}  \frac{\widehat{\psi'_+}(\xi_1+is)\varphi(\xi_1)}{i(\xi_1-is)}   =
 -i\pi\delta(\xi_1)+\PV \frac{\widehat{\psi_+'}(\xi_1)\varphi(\xi_1)}{\xi_1}. 
\end{equation}

To apply Isozaki's condition, fix any monotone function
$\chi_-\in\cC^{\infty}(\bbR)$ which vanishes for $t > 2\epsilon$ and equals $1$
for $t < \epsilon$, where $0<\epsilon \ll \xi^\circ$.  Now suppose that $x$ lies
in a small cone $C_\nu = \{ |x_2| \leq \nu \, x_1\}$ (or
equivalently, $\{|\omega_2| \leq \nu \omega_1\}$ where $\omega = x/|x|$).  Then
$\chi_-(\omega \cdot \xi)$ is supported in a union of half-planes $\omega \cdot
\xi \leq 2\epsilon$ which is a translate by $(2\epsilon, 0)$ of the conic
sector dual to $C_\nu$.  This intersects the line $\{\xi_1 = \xi^\circ\}$,
hence it is not obvious that the function
\begin{equation}
Au(x)=\int\chi_-(\omega \cdot \xi)\widehat{\psi_+u}(\xi)  e^{ix\cdot\xi}\, d\xi 
\label{111}
\end{equation}
is rapidly decaying for $x \in C_\nu$.

Reinserting the translation in $\xi_1$ by $\xi^\circ$, the $\delta$-function
contribution from~\eqref{eqn121.207} is
\begin{equation}\label{eqn98.202}
-i\pi e^{ix_1\xi^\circ}\int\chi_-\left(\omega_1\xi^\circ+\omega_2\xi_2\right)\hv(\xi_2)  e^{ix_2\xi_2}\, d\xi_2. 
\end{equation}
The condition $\omega \in C_\nu$ implies that $\omega_1>1/\sqrt{1+\nu^2}.$
Having fixed $\nu$, now choose $\epsilon$ so that $2\epsilon <
\xi^\circ/\sqrt{1+\nu^2}$.  If $\omega_2 > 0$ (the case $\omega_2 < 0$ is
handled similarly), then the integrand is supported in the region $S_\omega =
\{\xi_2: \xi_2 < (2\epsilon - \omega_1 \xi^\circ)/\omega_2\}$.  
As $\omega_2 \neq 0$, we can integrate by parts $N$ times, for any $N > 0$, and
see that \eqref{eqn98.202} equals
\begin{equation}
\frac{1}{(-ix_2)^N}  e^{ix_1\xi^\circ}\int_{S_{\omega}}\pa_{\xi_2}^N\left[\chi_-\left(\omega_1\xi^\circ+\omega_2\xi_2\right)\hv(\xi_2)\right]  e^{ix_2\xi_2}\, d\xi_2.
\end{equation}
As $x_2=r\omega_2,$ this, in turn, is bounded in absolute value by
\begin{equation}
\begin{split}
C_{M,N}r^{-N} \omega_2^{-N} & \int_{-\infty}^{(2\epsilon- \xi^\circ \omega_1)/\omega_2} (1 + |\xi_2|)^{-M}\, d\xi_2 \\
& \leq C_{M,N}' r^{-N} \omega_2^{M+1-N} \frac{1}{(\xi^\circ \omega_1 - 2\epsilon)^{M+1}} 
\end{split}
\end{equation}
for any $M \geq 0$.  Since $1/(\xi^\circ \omega_1 - 2\epsilon) $ is fixed, and
bounded away from $0$ if $|\omega_2|<\nu;$ choosing $M=2N,$ we see that this
term is $O((\omega_2/r)^N)$ as $r\to\infty.$

To estimate the contribution of the principal value term, set $t = \xi_1 - \xi^\circ$, so $|t| < 2\delta$ in the support of $\varphi$. 
Then we must examine
\begin{equation}\label{Bu}
\begin{split}
 Bu(x) & =e^{i\xi^\circ  x_1}\int\hv(\xi_2)e^{ix_2\xi_2}\\  & \left[\PV\int\limits_{-2\delta}^{2\delta}
\frac{\chi_-\left(\omega_1\xi^\circ+\omega_1 t +\omega_2\xi_2\right)\widehat{\psi_+'}(t)\varphi(t) e^{ix_1 t}}{t}\, dt\right]d\xi_2,
\end{split}
\end{equation}
We can now argue much as before: the term in square brackets is a distributional pairing
\[
\big\langle \PV \frac{1}{t}, \chi_-(\omega_1 \xi^\circ + \omega_1 t + \omega_2 \xi_2) \widehat{\psi'_+}(t) \varphi(t) e^{ix_1 t} \big\rangle.
\]
The function on the right of this pairing is smooth in $t$, uniformly
supported in $[-2\delta, 2\delta]$, and uniformly bounded as a function of
$\xi_2$.  Hence the pairing itself is smooth and uniformly bounded in $\xi_2$.
Finally, then, the integral in $\xi_2$ in \eqref{Bu} is the inverse Fourier
transform of a Schwartz function, hence it too is Schwartz, as claimed.  This
completes the proof that a wave-guide mode satisfies Isozaki's form of the outgoing condition.

As noted, the solutions constructed in \cite{EpWG2023_1, EpWG2023_2} also have
continuous spectral contributions, which are denoted by $u_0^{l,r},u_1^{l,r},
u_{c0}^{l,r}, u_{c1}^{l,r},$  where $\pm x_1>0,$
\begin{equation}
  \begin{split}
    u^{l,r}_0(x_1,x_2)=\cS_k\tau(x_1,x_2),\quad &u^{l,r}_1(x_1,x_2)=\cD_k\sigma(x_1,x_2),\\
  u^{l,r}_{c0}(x_1,x_2)=\cW^{l,r}_c\tau(x_1,x_2),\quad &u^{l,r}_{c1}(x_1,x_2)=\cW^{l,r\, '}_c\sigma(x_1,x_2).
  \end{split}
\end{equation}
See equations (95) and (208) in~ \cite{EpWG2023_2}.  \Rd {We must show that these too
satisfy the outgoing conditions. In order for this to be true, the densities,
$(\sigma,\tau),$ in~\eqref{eqn103.700} must have asymptotic expansions like
those in~\eqref{eqn114.220}. This condition is satisfied if the data $(g,h)$
in~\eqref{eqn88.220} have asymptotic expansions of this type as well. If
$u^{\In}_{l,r}$ are defined by wave-guide modes, point sources, or packets of
plane wave then this will be true. See Section 6 of~\cite{EpWG2023_1}.}\Bk

The estimates quoted below are the content of Theorem 2 from~\cite{EpWG2023_2},
and hold for $j = 0, 1$. First,
\begin{equation}
(\pa_{r}-ik)u_{j}^{l,r}(r\eta)=O(r^{-\frac 32}) \ \mbox{in the half-planes}\  \pm \eta_1 > 0,
\end{equation}
uniformly as $\eta_2\to \pm 1.$  Hence the $u_j^{l,r}$ are outgoing in the classical sense everywhere, 
including within the channels.  The  normal symbol $N_0(\pa_r-ik)=i(\tau-k)$ is non-zero on the characteristic set, except on $R^+_{k^2}.$  Thus away from $\eta_1=0,$ for some $m,$
\begin{equation}
\WF^{m,0}_{\Sc}(u_{j}^{l,r})\restrictedto_{\{\eta_1\neq 0\}}\subset R^+_{k^2,\mf}\cup \{\tau=k\}_{\ff}.
\end{equation}
Melrose's and Vasy's propagation  results then imply that
  \begin{equation}
\WF_{\Sc}(u_{j}^{l,r})\restrictedto_{\{\eta_1\neq 0\}}\subset R^+_{k^2,\mf}\cup \{\tau=k\}_{\ff}.
  \end{equation}
  See point 3. before Remark~\ref{rmk4.202}. 
We return to the behavior above $\eta_1=0$ below.

\begin{figure}[H]
\centering \includegraphics[width= 8cm]{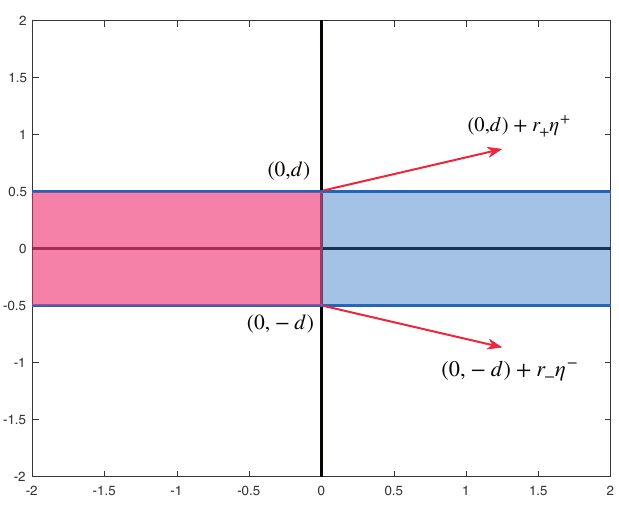}
\caption{Polar coordinates, $r_{\pm}\eta^{\pm},$  adapted to the locations of the channels. }
\label{adp_pol_coord}
\end{figure}

We now use polar coordinates $r_{\pm}=\sqrt{x_1^2+(x_2\mp d)^2}$, adapted to the
channel, which is assumed to lie in the set $\{(x_1,x_2):\: -d\leq x_2\leq d\},$
see Figure~\ref{adp_pol_coord}.  In terms of these coordinates there are
outgoing estimates
\begin{equation}
(\pa_{r_{\pm}}-ik)u_{cj}^{l,r}(r_{\pm}\eta^{\pm}\pm (0,d))=O(r_{\pm}^{-\frac 32})
\end{equation}
for $0<\pm\eta^{\pm}_{2}<1$; these are uniform when $\pm\eta^{\pm}_1\to 1^-$ and $\pm\eta^{\pm}_1\to 0^+$.   Moreover, 
when $x_2$ is bounded, 
\begin{equation}
(\pm\pa_{x_1}- ik)u_{cj}^{l,r}(x_1,x_2)=O(x_1^{-\frac 32}).
\end{equation}
These estimates are proved in Theorem 3, and Remark 12 of~\cite{EpWG2023_2}.

Choose a monotone function $\psi_{+}\in\cC^{\infty}(\bbR)$ with $\psi_+(x_2) = 0$ for $x_2 < d+1$ and $1$ for $x_2 > d+2$,
and set
\begin{equation}
      \psi_0(x_2)=1-(\psi_+(x_2)+\psi_+(-x_2)).
\end{equation}
The operators
\begin{equation}\label{eqn198.219}
D_{\pm}=\psi_+(x_2)\pa_{r_+}\pm\psi_0(x_2)\pa_{x_1}+\psi_+(-x_2)\pa_{r_-}\in\Psi_{3\Sc}^{1,0}(\overline{\bbR^2_{\pm}}_{\WG}),
\end{equation}
and have normal symbols $i\tau.$ From the estimates above it follows that
\begin{equation}
  (D_{\pm}-ik)u^{l,r}_{cj}(r\eta)=O(r^{-\frac 32})\text{ for }j=0,1.
\end{equation}
Using Vasy's description of the wave-front set, this implies for some $m\in\bbR,$
\begin{equation}
\WF^{m,0}_{3\Sc}(u_{cj}^{l,r})\restrictedto_{\bbR^2_{\pm}}\subset \left[R^+_{k^2,\mf}\cup \{\tau=k\}_{\ff}\right]_{\{\pm x_1>0\}}.
\end{equation}
Hence $u_{cj}^{l,r}\restrictedto_{\bbR^2_{\pm}}$ satisfy Vasy's outgoing radiation condition and therefore
\begin{equation}
\WF_{3\Sc}(u_{cj}^{l,r})\restrictedto_{\bbR^2_{\pm}}\subset \left[R^+_{k^2,\mf}\cup \{\tau=k\}_{\ff}\right]_{\{\pm x_1>0\}}.
\end{equation} 

To verify the outgoing conditions for these terms in the sense of
Isozaki, we refer to Appendix~\ref{App1}. There we compute the conic Fourier
transform of the leading term in the asymptotic expansion of a classically
outgoing solution and show that it satisfies Isozaki's outgoing condition. The
higher order terms of the expansion belong to $L^2(\bbR^2)$ and hence also
satisfy Isozaki's condition.
   
It remains then to discuss the wave-front set over the points where $\eta_1=0.$
As a brief preview of the argument, the key point, proved in
\cite{VasyAsterisque}, is that the scattering wave-front sets of formal solutions have
a certain regularity which is obtained from a propagation of singularities
theorem. This regularity shows that there cannot be any `isolated' parts of this
wave-front set over $\eta_1 = 0$. This is further explained in the proof of
Theorem~\ref{thm3}.

\subsection{Uniqueness for the Integral Equations}
The foregoing analysis  leads to the expected uniqueness result for the integral
equations introduced in~\cite{EpWG2023_1} for the solution of the transmission
problem in~\eqref{eqn88.220}. As noted above, the solutions to the transmission
problem are expressed as
\begin{equation}
u^{l,r}=-\cE^{l,r\, '}\sigma+\cE^{l,r}\tau,
\end{equation}
where
\begin{equation}
\cE^{l,r}\tau=\cS_{k}\tau+\cW^{l,r}\tau,\text{ and } \cE^{l,r\,'}\sigma=\cD_{k}\sigma+\cW^{l,r\, '}\sigma.
\end{equation}
The integral equations along $\{x_1=0\}$ then take the form 
\begin{equation}\label{eqn122.203}
     \left(\begin{matrix}\Id & D\\ C&\Id\end{matrix}\right)
       \left(\begin{matrix}\sigma\\ \tau\end{matrix}\right)=
         \left(\begin{matrix}g\\h\end{matrix}\right).
   \end{equation}
We refer to Section 5 of~\cite{EpWG2023_1} for further details. Using the uniqueness of the outgoing solution we can finally 
complete the analysis of this system of equations.

To state this result we need to define some Banach spaces.
\begin{definition}
 For $\gamma\in\bbR,$ let
  \begin{equation}
    \cC_{\gamma}(\bbR)=\{f\in\cC^0(\bbR):\:|f|_{\gamma}=\sup\{(1+|x|)^{\gamma}|f(x)|\}<\infty:\: x\in\bbR\}.
    \end{equation}
\end{definition}

\begin{theorem}\label{thm3}
For any $0<\gamma<\frac 12,$ and $(g,h)\in\cC_{\gamma}(\bbR)\oplus
\cC_{\gamma+\frac 12}(\bbR)$ the system of equations in~\eqref{eqn122.203} has a
unique solution $(\sigma,\tau)\in\cC_{\gamma}(\bbR)\oplus \cC_{\gamma+\frac
  12}(\bbR).$
\end{theorem}
\begin{remark}
The proof is quite similar to the proof of Theorem 3.41 in~\cite{ColtonKress}. 
\end{remark}
\begin{proof}
 In~\cite{EpWG2023_1} it is shown that the equations in~\eqref{eqn122.203}
     are Fredholm  equations of index zero on the spaces
     $\cC_{\gamma}(\bbR)\oplus \cC_{\gamma+\frac 12}(\bbR),$ for any
     $0<\gamma<\frac 12.$ It therefore suffices to show that the null-space is
     trivial. Let $(\sigma,\tau)$ satisfy~\eqref{eqn122.203} with $(g,h)=(0,0).$
     In this circumstance Theorem 1 of~\cite{EpWG2023_2} shows that $\sigma$ and
$\tau$ have asymptotic expansions as $\pm x_2\to\infty,$
\begin{equation}\label{eqn114.220}
         \begin{split}
   &\sigma(x_2)=\frac{e^{ik|x_2|}}{|x_2|^{\frac
       12}}\sum_{l=0}^N\frac{a^{\pm}_l}{|x_2|^l}+O\left(|x_2|^{-N-\frac
         32}\right),\\
       &\tau(x_2)=\frac{e^{ik|x_2|}}{|x_2|^{\frac
       32}}\sum_{l=0}^N\frac{b^{\pm}_l}{|x_2|^l}+O\left(|x_2|^{-N-\frac
     52}\right), 
     \end{split}
     \end{equation}
for any $N>0.$ The existence of these expansions then implies that the solutions $u^{l,r}$ satisfy the outgoing estimates 
in Theorems 2 and 3 of~\cite{EpWG2023_2}. Given the homogeneous transmission condition, the function
\begin{equation}
u^{\tot}(x_1,x_2)=
\begin{cases}
       u^{l}(x_1,x_2)\text{ for }x_1\leq 0,\\
         u^{r}(x_1,x_2)\text{ for }x_1\geq 0
\end{cases}
\end{equation}
defines a weak solution to $(\Delta+q(x_1,x_2)+k^2)u^{\tot}=0,$ smooth away from
$\{x_1=0,|x_2|<d\}$.

The analysis above shows that $\WF_{3\Sc}(u^{\tot})\restrictedto_{\{\eta_1\neq
  0\}}\subset R^+_{k^2,\mf}\cup \{\tau=k\}_{\ff}.$ Melrose's propagation results
exclude the possibility of a non-trivial intersection with
$\Sigma_{k^2,\mf}\setminus R^-_{k^2,\mf}\cup R^+_{k^2,\mf},$ but, a priori, it
is possible that $\WF_{3\Sc}(u^{\tot})$ over $\eta_1=0$ intersects
$R^-_{k^2,\mf}$ non-trivially.  Note that $q(x)=0$ in conic neighborhoods of
$\{\eta_2=\pm 1\}$ and therefore Proposition 2.8 of~\cite{VasyJFA97} shows that
$u^{\tot}$ has asymptotic expansions, in neighborhoods of $\eta_2=\pm 1,$ of the
form
\begin{equation}
  u^{\tot}(r\eta)\sim \frac{e^{ikr}}{r^{\frac
      12}}\sum_{j=0}^{\infty}\frac{a_j^+(\eta)}{r^j}+
  \frac{e^{-ikr}}{r^{\frac 12}}\sum_{j=0}^{\infty}\frac{a_j^-(\eta)}{r^j},
 \end{equation}
with $a^{\pm}_j(\eta)$ smooth functions of $\eta.$ Since all terms
$a^-_j(\eta)=0,$ for $\eta_2\neq \pm 1,$ the $e^{-ikr}$-term is absent, and
$u^{\tot}$ is everywhere outgoing and therefore $u^{\tot}\equiv 0.$ This is
  Proposition 17.8 in~\cite{VasyAsterisque}, which is restated above as
  Theorem~\ref{thm1}. 

To show that the data must be zero we interchange the roles of left and right, defining
            \begin{equation}
              \begin{split}
                v^r(x_1,x_2)&=-\cE^{l\,'}\sigma(x_1,x_2)+\cE^l\tau(x_1,x_2)\text{ for
                }x_1\geq 0,\\
                 v^l(x_1,x_2)&=\cE^{r\,'}\sigma(x_1,x_2)-\cE^r\tau(x_1,x_2)\text{ for
                }x_1\leq 0.
              \end{split}
            \end{equation}
  Note  the sign change in the definition of $v^l.$ As the functions
  $$\cW^{l,r}f(x_1,x_2),\, \cW^{l,r\,'}f(x_1,x_2),\,\pa_{x_1}\cW^{l,r\,'}f(x_1,x_2)$$
  are continuous across $x_1=0,$ it is elementary to show, using the classical jump
  relations for $\cS_{k}$ and $\cD_{k},$ that
  \begin{equation}\label{eqn127.203}
    v^r(0^+,x_2)=v^l(0^-,x_2)=-\sigma(x_2)\text{ and }\pa_{x_1}v^r(0^+,x_2)=\pa_{x_1}v^l(0^-,x_2)=\tau(x_2).
  \end{equation}
  The analysis used to show that $u^{\tot}$ is outgoing applies equally well to
     \begin{equation}
     v^{\tot}(x_1,x_2)=
     \begin{cases}
       v^{l}(x_1,x_2)\text{ for }x_1\leq 0,\\
         v^{r}(x_1,x_2)\text{ for }x_1\geq 0,
     \end{cases}
     \end{equation}
     which also satisfies  homogeneous jump conditions.  This function is
     therefore an outgoing solution to
     \begin{equation}
       (\Delta+q(-x_1,x_2)+k^2)v^{\tot}=0,
     \end{equation}
     which must also vanish. From~\eqref{eqn127.203} it follows that $\sigma=\tau=0.$
     \end{proof}

\subsection{Scattering via the Limiting Absorption Principle}\label{sec.lap_scat}
Solutions to scattering problems can also be found using the limiting absorption
  principle, (Theorem~\ref{thm2}), however the data in this formulation is used
  differently from how it used in the transmission formulation above. 
We can use an argument, similar to that used in the proof of Theorem~\ref{thm3},
to show that the solution $u^{\tot}$ defined in~\eqref{eqn90.220}, for certain
types  of incoming data, agrees with the limiting absorption
solution. The simplest case is when $u^{\In}_r=0$ and
\begin{equation}
  u^{\In}_l(x_1,x_2)=\sum_{j=0}^{N_l}c_lv_l(x_2)e^{ix_1\sqrt{k^2+E_{l,j}^2}},
\end{equation}
is a sum of incoming wave-guide modes from the left.  In this case the integral
equations in~\eqref{eqn122.203} are solvable and the solutions $(\sigma,\tau)$
have asymptotic expansions like those in~\eqref{eqn114.220}.  The functions
$u^{l,r}$ defined in~\eqref{eqn89.220}, which are solutions to
$(\Delta+q_{l,r}+k^2)u^{l,r}=0,$ are outgoing in their respective half spaces.

To employ the limiting absorption principle,  choose a function
$\varphi\in\cC^{\infty}(\bbR^2),$ supported in $x_1<0,$ which equals 1 in a
conic neighborhood $(-\infty,-1)\times \{0\}.$ With $w=(\Delta+k^2+q_l)[\varphi
  u^{\In}_l]\in\cS(\bbR^2),$ we let
\begin{equation}
  u^{\out}=-(\Delta+k^2+q+i0^+)^{-1}w
\end{equation}
be the unique outgoing solution, and  set $u_{\LAP}=u^{\out}+\varphi u^{\In}_l.$ It is
clear that
\begin{equation}
  (\Delta+k^2+q)(u_{\LAP}-u^{\tot})=0,
\end{equation}
and $u_{\LAP}-u^{\tot}$ is outgoing everywhere, except possibly over
$\{\eta_2=\pm 1\}.$ The argument used in the proof of Theorem~\ref{thm3} applies
to show that, in fact, this difference is everywhere outgoing, and therefore
$u_{\LAP}\equiv u^{\tot}.$

A similar argument applies with data defined by an incoming wave packet, as
described in Section 6 of Part I. As the description of this data would require
a lengthy discussion, we leave this case to the interested reader.  Generally
speaking, the solutions found using the integral equations in~\eqref{eqn90.220}
are outgoing provided the data $(g,h)$ is itself outgoing. This means that it
has a finite order asymptotic expansion, like (60) or (61) in~\cite{EpWG2023_2},
which implies that the sources $(\sigma,\tau)$ also have such expansions. While
we have not given the details, finite order expansions suffice because we only
need to know that~\eqref{eqn85.220} or, equivalently, \eqref{eqn108.213} holds
for an $l>-\frac 12.$

\section{Conclusion}
This paper concludes our introduction of mathematical foundations for the
problem of scattering scalar waves from open wave-guide networks. We have shown that this
problem has  physically motived, practically verifiable radiation conditions, which imply uniqueness,
and that the limiting absorption principle holds as
well. In~\cite{EpWG2023_1,EpWG2023_2} we have presented an effective method to
solve a particular model wave-guide network problem, which has been implemented
numerically in~\cite{GE_2024}. 

There are obviously many directions for further developments in this field. An
important applied problem is the development of more flexible numerical methods
that apply to more complex geometries in 2 and 3 dimensions. Once numerically
implemented, the method presented in~\cite{EpWG2023_1,EpWG2023_2} can serve as a
`gold standard' for new methods that may be more difficult to  rigorously
analyze. We also hope to establish the existence of  complete asymptotic expansions for outgoing
solutions, which are valid in neighborhoods of the channel ends. 

Another important direction is the study of the full scattering operator defined
by the wave guide problem. While its general properties can be obtained using
the geometric microlocal methods of Melrose and Vasy, it will be necessary, in
applications, to develop more explicit methods to estimate the relative sizes of
the different components of the scattering operator, i.e., to measure how much energy remains
in the channels vs.\  radiation escaping into the free regions. This could be quite
useful for the practical design of `optimal' wave-guide networks. In
Appendix~\ref{chan-to-chan} we show that the channel-to-channel scattering coefficients are
well defined.

The methods introduced here readily generalize to study the open wave guide
problem for Maxwell's equations.  On the free boundary, the scattering
wave-front set readily generalizes to this setting, as does the notion of
polarization of the principal singularities, see~\cite{Dencker1982}. The classical
Silver-M\"uller radiation conditions are easily understood in this language. It
seems quite likely that channels can again be accommodated, though the details
have not yet been worked out.

\appendix
\section{Appendix: Conic Fourier Transform of a Classically Outgoing Solution}\label{App1} 
In this appendix we carry out a relatively elementary but lengthy computation
showing that the leading term in the asymptotic expansion of a classically
outgoing solution in two dimensions satisfies Isozaki's radiation
condition. In~\cite[Lemma 1.4]{Isozaki94} Isozaki proves this abstractly.
The corresponding fact is also straightforward in Vasy's formulation using the
observation that $\WF_{3\Sc}(v)\subset\Sigma_{k^2}\cap \{\tau=k\}=R^+_{k^2}.$

This leading term takes the form $v(r,\theta) = a(\theta)e^{ik r}/\sqrt{r},$
where $a \in \cC^\infty(\pa\overline{\bbR^2})$.  We compute its conic Fourier
transform around any ray $r \zeta(\theta_0)$, $r > 0$, where  we let
$$
\zeta(\theta)=(\cos\theta,\sin\theta).
$$
For simplicity of notation we may as well assume that $\theta_0 = 0$, which we do below.

Multiply $v$ by smooth localizing functions $\psi_+(r)$ and $\varphi(\theta)$;
here $\psi_+(r)$ vanishes for $r < 1$ and equals $1$ for $r > 2$, and $\varphi$
is supported in $(-\delta,\delta)$ with $\varphi(\theta)=1$ when
$|\theta|<\delta/2.$ The integral defining the conic Fourier transform converges if we
replace $k$ by $k + is$, $s > 0$, to get
\begin{equation}
\hv_{s}(\xi)=
\int\frac{a(\theta)\psi_+(r)\varphi(\theta)e^{i(kr-y\cdot\xi)}e^{-rs}}{\sqrt{r}}\,
r dr d\theta.
\end{equation}
Our goal is to compute the distributional limit $\hv=\lim_{s\to
  0^+}\hv_{s}.$  We give the following description of this limit:
\begin{lemma} Given $\nu>0,$ there exists a $\delta>0$ so that, with
  $\supp\varphi\subset (-\delta,\delta),$ the distributional limit $\hv$ is
  smooth in the complement of $B_\nu( k,0)$, and is rapidly decreasing, along
  with all derivatives, as $|\xi| \to \infty.$
\end{lemma}

Denoting by $\chi_-$ the same function as in \eqref{111}, using this  lemma we
show below that, with $0<\epsilon$ and $0<\delta$ chosen small enough,
\begin{equation}\label{eqn161.208}
Av(x)=  \int\chi_-\left(\frac{x}{|x|}\cdot\xi\right) \hv(\xi)e^{ix\cdot\xi}\, d\xi
\end{equation}
decreases rapidly for $x$ in a cone $C_{\delta'}=\{r\zeta(\theta):\:
|\theta|<\delta',\ r>0\}$. This is Isozaki's outgoing radiation condition.

\begin{proof}[Proof of the lemma]
To prove the lemma we set $x=r\zeta(\theta)$ and $\xi=\rho\zeta(\phi)$, so that
\begin{equation}
\hv_{s}(\xi)=  \int_{-\delta}^{\delta} \int_{1}^{\infty}a(\theta)\psi_+(r)\varphi(\theta)e^{ir(k-\rho\cos(\theta-\phi))}e^{-rs}\sqrt{r}\, drd\theta.
\end{equation}
Using $\pa_{\theta}e^{-ir\rho\cos(\theta-\phi)}=ir\rho\sin(\theta-\phi)e^{-ir\rho\cos(\theta-\phi)}$, 
we can integrate by parts in $\theta$ 
as often as we like to obtain
\begin{equation}
\begin{split}
\hv_{s}(\rho,\phi) = \int_{1}^{\infty}   \int_{-\delta}^{\delta} & \left(\pa_{\theta}\frac{-1}{ir\rho\sin(\theta-\phi)}\right)^{N}
(a(\theta)\varphi(\theta))\\ & e^{ir(k-\rho\cos(\theta-\phi))}\psi_+(r)e^{-rs}\sqrt{r} \, d\theta\,dr.
\end{split}
\end{equation}
Now let $s\to 0^+$ to conclude that $\hv(\rho,\phi)$ is rapidly decreasing and smooth in the double cone
$$
V_{\delta}=\{(\rho,\phi):\:|\phi|>2\delta\ \mbox{and}\ |\phi-\pi| > 2\delta,\ \rho\neq 0\}. 
$$

Next, since 
\[
\pa_r e^{ir(k-\rho\cos(\theta-\phi) + is)}=i(k-\rho\cos(\theta-\phi) + is )e^{ir(k-\rho\cos(\theta-\phi) + is)},
\]
we can integrate by parts in $r$ as often as we like when $\rho<k$ to get
\begin{equation}
\begin{split}
\hv_{s}(\rho,\phi) & = \int_{-\delta}^{\delta} \int_{1}^{\infty}(a(\theta)\varphi(\theta)) e^{ir(k-\rho\cos(\theta-\phi))} e^{-rs} \\ & 
\left[\frac{-1}{i(k-\rho\cos(\theta-\phi) + is)}\,\pa_r\right]^N (\psi_+(r)\sqrt{r})\, dr d\theta.
\end{split}
\end{equation}
The resulting integrand is bounded independently of $s\geq 0$ by $C_Nr^{\frac 12-N}$, thus we can let $s\to 0^+$
and again conclude that $\hv(\rho,\phi) \in \cC^\infty$ when $\rho<k$.

Next, if $|\phi-\pi|<2\delta$, then $\cos(\theta-\phi)<\cos(3\delta-\pi)$, and this is strictly negative when $\delta\ll 1$,
so we can again integrate by parts in $r$ to conclude that $\hv(\rho,\phi)$ is smooth and rapidly decreasing in the cone 
$\{|\phi-\pi|\leq 2\delta\}.$
Finally, when $|\phi|\leq 2\delta,$ then $\rho\cos(\theta-\phi)>\rho\cos(3\delta).$ If $\rho>\frac{k}{\cos(3\delta)},$ then a
final integration by parts shows that $v(\rho,\phi)$ is smooth and rapidly decreasing in the set $\{ |\phi|\leq 2\delta, 
\rho>\frac{k}{\cos(3\delta)}\}.$ 

Combining all these estimates, we conclude that for any $\nu>0,$ there exists a $\delta>0$ so that 
$\hv\in\cC^{\infty}(\bbR^2\setminus B_{\nu}(k,0) ),$ and is rapidly decreasing with all derivatives 
as $\rho\to\infty.$
\end{proof}

Now consider $Av(x)$ as in~\eqref{eqn161.208}, where $x\in C_{\delta'}=\{r\zeta(\theta):\: |\theta|<\delta',r>0\}$
for some appropriate $\delta'>0.$  The integrand is supported in $\{\xi:\omega\cdot\xi<2\epsilon\}.$ 
Thus by choosing $\nu\ll k,$ and $\epsilon<\nu,$ we can then choose $\delta'>0$ so that when $x\in
C_{\delta'},$ the integrand in $Av(x)$ is a smooth, rapidly decreasing function of $\xi$ with support a fixed 
positive distance from $B_{\nu}(k,0).$  Thus $Av(x)$ is also rapidly decreasing as $|x|\to\infty$, so $v$ 
satisfies Isozaki's outgoing condition.

Referring back to the notation of Section~\ref{sec_2d}, these calculations suffice to handle the terms arising
from $u_{0}^{l,r}, u_{1}^{l,r}$ and $u_{c0}^{l,r}, u_{c1}^{l,r}$ away from the ends of the channels. Similar, though 
more complicated, calculations handle the contributions involving $u_{c0}^{l,r}, u_{c1}^{l,r}$ in neighborhoods of 
the channels. We leave these estimates to the interested reader.

\section{Channel-to-Channel Scattering Coefficients}\label{chan-to-chan}
A fundamental question of scattering theory for open wave-guide networks is to
understand how an incoming wave-guide mode is scattered by the network of
wave-guides into a sum of outgoing wave-guide modes and radiation. Here we
consider the wave-guide mode portion of the scattered field. To formulate this
question mathematically, we now construct the appropriate class of solutions
using the various tools we have been discussing.  For each $\alpha\in\cA,$
choose a conic cutoff function $\varphi_\alpha$ which equals $1$ in the exterior
of a large ball $B_{2R}(0)$ intersected with a conic region around $v_\alpha,$
which has support outside $B_R(0)$ and conic neighborhoods of all the other
$v_\beta$.  Thus $q(x) = q_\alpha(x^\alpha)$ in $\supp\varphi_{\alpha}.$

Now consider an incoming wave-guide mode associated to this channel and
localized to this neighborhood; applying the operator to it yields a rapidly
vanishing function
\begin{equation}
f_{\alpha,j}=(\Delta+q+k^2)\left[\varphi_{\alpha}(x)e^{-ix_{\alpha}\sqrt{E^2_{\alpha,j}+k^2}}u_{\alpha,j}(x^{\alpha})\right] \in \cS(\bbR^d).
\end{equation}
We now apply the outgoing resolvent to obtain the outgoing solution,
$v^+_{\alpha,j},$ to
\begin{equation}\label{eqn98.206}
(\Delta+q+k^2)v^+_{\alpha,j}=f_{\alpha,j}.
\end{equation}
The difference $w_{\alpha,j}=e^{-ix_{\alpha}\sqrt{E^2_{\alpha,j}+k^2}}u_{\alpha,j}(x^{\alpha})-v^+_{\alpha,j}$ is then a generalized eigenfunction: 
\begin{equation}
(\Delta+q+k^2)w_{\alpha,j}=0.
\end{equation}
It has an incoming component along the $v_\alpha$ channel and a superposition of
outgoing components in all of the channels, along with outgoing radiation.

We seek a formula for the outgoing guided-mode contributions.  We can isolate
the component corresponding to any wave-guide mode $u_{\beta,l}(x^{\beta})$ with
energy $E^2_{\beta,l}$ by taking the inner product localized along that channel. This
gives a function
\begin{equation}
s_{\alpha,j;\beta,l}(x_{\beta})=\int_{\bbR^{d-1}}\varphi_{\beta}(x_{\beta},x^{\beta}) v^+_{\alpha,j}(x_{\beta},x^{\beta})\overline{u_{\beta,l}(x^{\beta})}dx^{\beta}
\end{equation}
which vanishes by construction when $x_\beta < R.$ These functions have a very
special form.

\begin{lemma}
  The function $s_{\alpha,j;\beta,l}(x_{\beta})$ is given by
  \begin{equation}\label{eqn101.205}
 s_{\alpha,j;\beta,l}(x_{\beta})=S_{\alpha,j;\beta,l}e^{ix_{\beta}\sqrt{E_{\beta,l}^2+k^2}}+\sigma(x_{\beta}), \qquad  \sigma \in \cS(\bbR).
\end{equation}
\end{lemma}

We can thus define the scattering coefficient $S_{\alpha,j;\beta,l}$ between the incoming mode $(\alpha,j)$ and outgoing mode $(\beta,l)$ by
\begin{equation}
S_{\alpha,j;\beta,l}=\lim_{x_{\beta}\to\infty}e^{-ix_{\beta}\sqrt{E_{\beta,l}^2+k^2}}s_{\alpha,j;\beta,l}(x_{\beta}).    
\end{equation}
We must verify~\eqref{eqn101.205} and show that this limit is independent of all choices made in its definition.

The independence of choice of cut-offs is easily established.  If $\{\varphi_{\alpha}':\:\alpha\in\cA\}$ is another collection of conic cutoff
functions, with corresponding outgoing solutions $v_{\alpha,j}^{+'}$ and projections $\{s_{\alpha,j;\beta,l}'(x_{\beta})\}$ as in~\eqref{eqn101.205},
then the difference $(v_{\alpha,j}^+-v_{\alpha,j}^{+'})$ is the unique outgoing solution to
\begin{equation}
(\Delta+q+k^2)v=(\Delta+q+k^2)\left[(\varphi_{\alpha}-\varphi_{\alpha}') e^{-ix_{\alpha}\sqrt{E^2_{\alpha,j}+k^2}}u_{\alpha,j}(x^{\alpha})\right].
\end{equation}
The difference on the right is rapidly decreasing, hence a fortiori outgoing. Thus by uniqueness of outgoing solutions, we see that
\begin{equation}
(v_{\alpha,j}^+-v_{\alpha,j}^{+'})= \left[(\varphi_{\alpha}-\varphi_{\alpha}')  e^{-ix_{\alpha}\sqrt{E^2_{\alpha,j}+k^2}}u_{\alpha,j}(x^{\alpha})\right],
\end{equation}
hence
\begin{equation}
s_{\alpha,j;\beta,l}(x_{\beta})-s_{\alpha,j;\beta,l}'(x_{\beta})\in\cS(\bbR), 
\end{equation}
so the coefficients $\{S_{\alpha,j;\beta,l}\}$ are well defined.
\begin{proof}[Proof of the lemma]
We next establish~\eqref{eqn101.205}. Using $(\Delta_{\beta}+q_{\beta}(x^{\beta}))u_{\beta,l}=E_{\beta,l}^2u_{\beta,l}$, we write
\begin{equation}
\begin{split}
s_{\alpha,j;\beta,l}(x_{\beta}) = & \\
\frac{1}{E_{\beta,l}^2} &\int\varphi_{\beta}(x_{\beta},x^{\beta})v^+_{\alpha,j}(x_{\beta},x^{\beta})
\overline{(\Delta_{\beta}+q_{\beta}(x^{\beta}))u_{\beta,l}(x^{\beta})}dx^{\beta}.
\end{split}
\end{equation}
The function
$\varphi_{\beta}(x_{\beta},x^{\beta})v_{\alpha,j}(x_{\beta},x^{\beta})$ is
smooth with bounded derivatives, and $u_{\beta,l}(x^{\beta})$ decays
exponentially, so we can integrate by parts to obtain:
\begin{equation}
s_{\alpha,j;\beta,l}(x_{\beta})=\frac{1}{E_{\beta,l}^2}\int(\Delta_{\beta}+q_{\beta}(x^{\beta}))\left[\varphi_{\beta}(x_{\beta},x^{\beta})v^+_{\alpha,j}(x_{\beta},x^{\beta})\right]
\overline{u_{\beta,l}(x^{\beta})}dx^{\beta}.  
\end{equation}
Using~\eqref{eqn98.206}, we see that
  \begin{multline}
    s_{\alpha,j;\beta,l}(x_{\beta})=-\frac{1}{E_{\beta,l}^2}(\pa_{x_{\beta}}^2+k^2)\int \varphi_{\beta} v^+_{\alpha,j} \overline{u_{\beta,l}(x^{\beta})}dx^{\beta}+\\
    \frac{1}{E_{\beta,l}^2}\int\Bigg[ \varphi_{\beta}f_{\alpha,j}+ \Delta\varphi_{\beta}v^+_{\alpha,j}+2\nabla\varphi_{\beta}\cdot\nabla
      v^+_{\alpha,j}-\varphi_{\beta}(q-q_{\beta})v^+_{\alpha,j}
      \Bigg]    \overline{u_{\beta,l}(x^{\beta})}dx^{\beta}.  
  \end{multline}
From this and the support properties of the $q_{\alpha}$ and conic cut-offs $\varphi_{\alpha}$, it follows that
\begin{equation}
(\pa_{x_{\beta}}^2+k^2+E_{\beta,l}^2) s_{\alpha,j;\beta,l}= g\in\cS(\bbR).
\end{equation}
A standard analysis of this one-dimensional problem shows that there are constants $a,b,$ and smooth functions $h_0, h_1$ supported
in the positive half-line with $h_0(t) = 1$ for sufficiently large $t$ and $h_1$ Schwartz class, such that
\begin{equation}
s_{\alpha,j;\beta,l}(x_{\beta})=h_0(x_{\beta})\left[ae^{ix_{\beta}\sqrt{E_{\beta,l}^2+k^2}}+be^{-ix_{\beta}\sqrt{E_{\beta,l}^2+k^2}}\right]+ h_1(x_{\beta}).
\end{equation}

To establish~\eqref{eqn101.205}, it remains to show that $b=0.$  For any integrable function $m(x^\beta)$, denote by
$\tm(\xi^{\beta})$ its $(d-1)$-dimensional Fourier transform. Since $u_{\beta,l}\in\cS(\bbR^{d-1})$ we can apply the Plancherel formula to rewrite
\begin{equation}
s_{\alpha,j;\beta,l}(x_{\beta})=\frac{1}{(2\pi)^{d-1}}\int_{\bbR^{d-1}}\widetilde{\varphi_{\beta}
v^+_{\alpha,j}}(x_{\beta},\xi^{\beta})\overline{\tu_{\beta,l}(\xi^{\beta})}d\xi^{\beta}.
\end{equation}
Taking the Fourier transform in $x_{\beta}$ now gives 
\begin{equation}\label{eqn163.221}
\begin{split}
\frac{1}{(2\pi)^{d-1}}\int_{\bbR^{d-1}} & \widehat{\varphi_{\beta} v^+_{\alpha,j}}(\xi_{\beta},\xi^{\beta})\overline{\tu_{\beta,l}(\xi^{\beta})}d\xi^{\beta}=\\
& a\cF(h_0e^{ix_{\beta}\sqrt{E_{\beta,l}^2+k^2}})+b\cF(h_0e^{-ix_{\beta}\sqrt{E_{\beta,l}^2+k^2}}) +\hh_1(\xi_{\beta}).
\end{split}
\end{equation}

It follows, as above in the analysis of~\eqref{eqn105.204}, that there are non-zero constants $c_0,c_1$  and a function $h_2 \in\cS(\bbR)$ so that
\begin{multline}\label{eqn114.207}
\cF(h_0e^{-ix_{\beta}\sqrt{E_{\beta,l}^2+k^2}})=\\c_0\delta\left(\xi_{\beta}+\sqrt{E_{\beta,l}^2+k^2}\right)+ 
c_1\PV\frac{\widehat{h_0'}(\xi_{\beta})}{\xi_{\beta}+\sqrt{E_{\beta,l}^2+k^2}}+h_2(\xi_{\beta}).
\end{multline}

Let $\chi_{-}\in\cC^{\infty}(\bbR)$ have its support in
$\left((-(k^2+E^2_{\beta,N'_{\beta}})^{\frac 12}-2\epsilon,-k+2\epsilon)\right),$
and be equal to 1 in $\left((-(k^2+E^2_{\beta,N'_{\beta}})^{\frac 12}
-\epsilon,-k+\epsilon)\right),$ for an $0<\epsilon\ll k/2.$ Using the outgoing
condition we see below that if $\varphi_{\beta}$ has sufficiently small conic support,
then
\begin{equation}\label{eqn161.205}
Pv^+_{\alpha,j}(x_{\beta},0)=\frac{1}{(2\pi)^{d-1}}\int_{\bbR}\int_{\bbR^{d-1}}\widehat{\varphi_{\beta}
 v^+_{\alpha,j}}(\xi_{\beta},\xi^{\beta})\chi_{-}(\xi_{\beta})\overline{\tu_{\beta,l}(\xi^{\beta})} e^{ix_{\beta}\xi_{\beta}}d\xi^{\beta}d\xi_{\beta}
\end{equation}
must tend to zero as $x_{\beta}\to\infty.$ To prove this statement requires
material from~\cite{VasyAsterisque} that we have not covered in detail. In the following
discussion we carefully follow the notation from this paper, but leave the
detailed verification to the reader.

The function $\overline{\tu_{\beta,l}(\xi^{\beta})}\in\cS(\bbR^{d-1});$
as described in Remark~\ref{rmk7.205}, we are free to replace the symbol of $P,$
$$p(x,\xi,x')=\varphi_{\beta}(x)\chi_-(\omega\cdot
\xi)\overline{\tu_{\beta,l}(\xi^{\beta})}\varphi_{\beta}(x'),$$ with
$p(x,\xi,x')/(1+|\xi|^2)^N,$ for any $N.$ Given $p\geq 0,$ by choosing $N$ large
enough we can arrange to have $p(x,\xi,x')/(1+|\xi|^2)^N$ satisfy any finite
number of symbolic estimates for an operator in
$\Psi^{-p,0}_{3\Sc}(\overline{\bbR^d}).$ Denote the associated operator by
$P_N.$

As $v^+_{\alpha,j}$ is outgoing, the $WF_{3\Sc}(v^+_{\alpha,j})\subset \{\tau>0\}.$
Moreover $WF_{3\Sc}(\Delta^r v^+_{\alpha,j})\subset WF_{3\Sc}(v^+_{\alpha,j}),$
for any $r\in\bbN,$ and therefore 
$$\WF'_{3\Sc}(P_N)\cap WF_{3\Sc}((\Id-\Delta)^N v^+_{\alpha,j})=\emptyset,$$
where $\WF'_{3\Sc}$ is the operator wave front set, see Definition 9.1
in~\cite{VasyAsterisque}. Arguing as in the proof of Lemma 9.8
from~\cite{VasyAsterisque}, we can show that for any $q,l\in\bbR$ the function
$$P_N(\Id-\Delta)^Nv^+_{\alpha,j}=Pv^+_{\alpha,j}\in
H^{q,l}(\overline{\bbR^d}).$$
Taking $q,l$ sufficiently large, we see that $P
v^+_{\alpha,j}(x)$ must tend pointwise to zero as $x\to\infty.$ The integral
in~\eqref{eqn161.205} is $P v^+_{\alpha,j}(x_{\beta},0),$ and therefore~\eqref{eqn163.221} and~\eqref{eqn114.207} imply that
\begin{equation}
  P v^+_{\alpha,j}(x_{\beta},0)=bc_0e^{-ix_{\beta}\sqrt{E_{\beta,l}^2+k^2}}+o(1).
\end{equation}
Hence the outgoing condition implies that $b=0.$
This completes the proof that the incoming-to-outgoing wave-guide mode
scattering coefficients, $\{S_{\alpha,j;\beta,l}:\:\alpha,\beta\in\cA \},$ are
well-defined. 
\end{proof}

These coefficients depend on the choice of origin for $\bbR^d$ as
well as the choice of orthonormal basis for the incoming and outgoing wave-guide
modes. Translating the origin by $x_0,$ has the effect of replacing
$S_{\alpha,j;\beta,l}$ by
\begin{equation}
e^{i\left[x_0\cdot v_{\alpha}\sqrt{E_{\alpha,j}^2+k^2}-x_0\cdot  v_{\beta}\sqrt{E_{\beta,l}^2+k^2}\right]}S_{\alpha,j;\beta,l},
\end{equation}
which is conjugation by a diagonal unitary matrix.  It remains a very interesting problem to give effective estimates for these
coefficients, and to bound the portion of the energy of the incoming signal that is dissipated in radiation.

The `full' scattering operator involves both radiation and wave-guide
modes. In~\cite{VasyAsterisque} the free-to-free part of the scattering operator
is described in considerable detail. In~\cite{Isozaki92} various other parts of
the scattering operator are described for the $N$-body Schr\"odinger case.
However, it is an interesting and important challenge to assemble all of this
into a more coherent and complete study of the full scattering operator for open
wave-guide networks.

\end{document}